%% file: statFEMDyn.tex
\documentclass[3p, number,sort&compress,times,fleqn, final]{elsarticle}

\makeatletter
\def\ps@pprintTitle{%
  \let\@oddhead\@empty
  \let\@evenhead\@empty
  \let\@oddfoot\@empty
  \let\@evenfoot\@empty
}
\makeatother

\graphicspath{{./figs/}}
\usepackage{url}
\usepackage{booktabs}
\usepackage{multirow}
\usepackage[table]{xcolor}
\usepackage{style/csml}
\usepackage{soul}
\usepackage{algorithm}
\usepackage{algpseudocode}

\usepackage[colorlinks = true,
linkcolor = blue,
urlcolor  = blue,
citecolor = blue,
anchorcolor = blue]{hyperref}

\newcommand{\red}[1]{\textcolor{red}{#1}}

\begin{document}
	
\begin{frontmatter}

\title{Statistical finite elements for sequential data synthesis in solid dynamics}

\author{Igor Kavrakov\fnref{fn1}}
\author{Yaswanth Sai Jetti\fnref{fn1}\corref{cor1}}
\ead{ysj22@cam.ac.uk}
\author{Ahmet Oguzhan Yuksel}
\author{Fehmi Cirak}
%\author{Fehmi Cirak\corref{cor1}}
\fntext[fn1]{These authors contributed equally to this work.}

\cortext[cor1]{Corresponding author.}

\address{Department of Engineering, University of Cambridge, Cambridge, CB2 1PZ, UK }

\begin{abstract}
We present an approach for synthesising observational data with elastodynamic finite element models by extending the statistical finite element method (statFEM) framework. The proposed formulation adopts a Bayesian filtering approach to account for uncertainties in the data, the finite element model, and the discrepancies between the model and the physical system. Observational data are assimilated while the state of the spatially discretised finite element problem is advanced in time using the stochastic variant of the explicit Newmark scheme. The prior probability density of the state is obtained by solving an incremental probabilistic forward problem, and the corresponding posterior density is obtained by conditioning on the data available at each time step. In the probabilistic forward problem, spatio-temporal Gaussian random fields representing the forcing, model misspecification, and material parameters are specified via their stochastic PDE formulation. The resulting non-Gaussian prior probability density is approximated using a perturbation approach, yielding a Gaussian posterior with closed-form mean and covariance. The hyperparameters of the random field representing model misspecification are calibrated by maximising the marginal likelihood of the data. The proposed approach is illustrated on one- and two-dimensional elastodynamic examples with synthetic data.

\end{abstract}
	

\end{frontmatter}

%\tableofcontents
%\newpage 
    
\input{introduction}
\input{forward}

\input{filtering}

\input{examples}

\input{conclusions}
\input{acknowledgements}

\appendix
\input{appendix}

\bibliographystyle{elsarticle-num-names}
\bibliography{statFEMDyn}

\end{document}

%% file: introduction.tex
%
%--------------------------------------------------------------------------------          
\section{Introduction \label{sec:introduction}}
%--------------------------------------------------------------------------------
%
The continuous monitoring of engineering systems and structures in operation requires computational techniques that integrate observational data with physics-based models~\cite{willcox2024role,wagg2020digital,sun2020review}. In practice, measurements are sparse and noisy, and computational models are subject to significant uncertainties due to modelling simplifications and assumptions. The conventional deterministic finite element method is ill-equipped to systematically incorporate (assimilate) observational data while accounting for uncertainties in a rational manner. The recently introduced statistical finite element method (statFEM) provides a probabilistic framework for incorporating data and modelling uncertainties in conventional finite element models; see~\cite{2021statistigirolamical} for the original formulation and~\cite{duffin2021statistical,febrianto2022digital,akyildiz2022statistical,koh2023stochastic,glyn2025statistical,narouie2025mechanical,karvonen2025error,narouie2026unsupervised} for subsequent developments. In statFEM, the uncertainties in the finite element model and the measurements are modelled using random fields and variables, with prescribed prior densities. Bayes' rule is used to determine the respective posterior densities, given the likelihood of the observed data under an assumed statistical observation model. A key feature of statFEM is the consideration of model misspecification to account for the inevitable discrepancies between the finite element model and the underlying physical system. Such discrepancies typically arise from misspecified initial and boundary conditions, the size of the computational domain, and uncertainties in material parameters or forcing.

In this work, we extend statFEM to elastic solid dynamic problems, providing a sequential Bayesian filtering formulation for data assimilation. Bayesian filtering addresses the inference of a spatio-temporal state from a limited set of indirect and noisy sequential observations~\cite{kaipio2006statistical,bishop2006pattern,bach2024machine}.  Observational data are assimilated while the solution of the spatially discretised finite element problem is advanced in time using a time-stepping procedure. The semi-discrete elastodynamic problem is expressed in state-space form so that the evolution of the state depends only on its value at the preceding time step (Markov property)~\cite{beck2010bayesian,sarkka2023bayesian}.  The prior probability distribution of the state is obtained by solving the incremental forward problem. This prior distribution is then updated by conditioning on the available observations. These two stages are commonly referred to as the prediction and update steps.  The mean of the resulting posterior distribution provides a point estimate of the inferred state, while credible intervals quantify the range of states consistent with the data and the prior.

The classical Kalman filter is a specific instance of Bayesian filtering for linear problems with Gaussian probability densities \cite{kalman1960new,jazwinski2007stochastic}. Extensions to nonlinear problems include the extended Kalman filter, based on first-order linearisation, and methods based on deterministic or Monte Carlo approximations, such as the unscented Kalman filter, ensemble Kalman filter, and particle filters~\cite{wan2000unscented,evensen1994sequential,doucet2001introduction}. Although Kalman filtering is primarily intended for estimating the state of a system, there are extensions for estimating the parameters or inputs of a system, such as the forcing or scalar constitutive parameters; see~\cite{corigliano2004parameter,lourens2012augmented,azam2015dual,maes2016joint,ebrahimian2018bayesian,caglio2025bayesian}.  These methods are typically applied to relatively low-dimensional dynamical systems or skeletal structures, such as frames, trusses, and coupled oscillators. The connection between Kalman filtering and statFEM has been explored by Duffin et al.~\cite{duffin2021statistical} for prototypical fluid dynamics problems, posed on relatively simple domains. Kalman filtering formulations for finite element models on non-trivial geometries with spatially distributed uncertainties remain to be developed. 

The statFEM construction naturally extends to sequential inference for transient elastodynamic problems.  In the proposed formulation, we discretise the elastodynamic forward problem in space using standard finite elements and integrate in time using the stochastic variant of the central difference, i.e. the explicit Newmark scheme~\cite{sarkka2019applied,burrageNumericalMethodsSecondOrder2007,mannellaQuasisymplecticIntegratorsStochastic2004a}.  Model misspecification and uncertainties in forcing and material parameters are parameterised as Gaussian spatio-temporal random fields. The forcing is decomposed into a known component and a residual component accounting for model misspecification, whose importance is well-established~\cite{kennedy2001bayesian,jiang2020sequential,novoa2024inferring}.  The spatial correlation structure of the random fields is specified via their stochastic PDE formulation~\cite{lindgren2011explicit,zhang2021stochastic,lindgren2022spde,koh2023stochastic}.  The respective stochastic PDEs are driven by spatial white noise, and are discretised on the same finite element mesh as the elastodynamic forward problem. The length scale, smoothness, and magnitude of the random fields are governed by the parameters of the stochastic PDEs.
The solution to the elastodynamic forward problem depends nonlinearly on the spatially correlated random material field, yielding a non-Gaussian probability density for the state in the prediction step. Although this density can be approximated using quadrature or sampling-based methods, such approaches are typically limited to relatively small problems~\cite{ghanem2003stochastic,xiu2010numerical,sudret2000stochastic,matthies2005galerkin}. We therefore approximate the probabilistic solution of the forward problem using a first-order perturbation approach \cite{liu1986random,sarkkaContinuousDiscreteCubatureKalman2012}, yielding a Gaussian approximation of the predictive density. The Gaussian approximation serves as the prior probability density, which is subsequently conditioned on the noisy observations in the update step. Owing to the posited observation model with additive noise, Bayes' rule yields a posterior probability density that remains Gaussian and is analytically tractable.  The hyperparameters of the stochastic random field representing model misspecification are calibrated by maximising the marginal likelihood (evidence) of the data. 
 
The remainder of the paper is organised as follows. Section \ref{sec:forward} develops the probabilistic forward model for elastodynamics, including the SPDE-based representation of material and forcing fields, spatial discretisation and perturbation-based propagation of mean and covariance. Section \ref{sec:bayesfil} introduces the Bayesian filtering framework for sequential data assimilation, derives the prediction and filtering equations for the augmented state formulation, and presents the empirical Bayes strategy for hyperparameter estimation. Section \ref{sec:examples} demonstrates the proposed methodology on representative one- and two-dimensional elastodynamic examples, highlighting state estimation, material field recovery and hyperparameter identification. Finally, Section \ref{sec:conclusions} summarises the main findings and discusses directions for future research.

%% file: forward.tex
\section{Probabilistic Forward Model for Elastodynamics \label{sec:forward}}
%--------------------------------------------------------------------------------
%
%In this section, we review the forward problem for elastodynamics formulated as a stochastic partial differential equation (SPDE) with random material properties and generalised forcing, accounting for applied loads and possible model discrepancies. The generalised forcing is modelled as a separable spatio-temporal random field, and the material property as a purely spatial random field. Both random fields are generated using auxiliary SPDEs that are distinct from the elastodynamic forward problem. We discretise in space using standard finite elements and in time using an explicit scheme. The mean and covariance of the solution field are obtained via first-order perturbation theory. 
In this section, we review the forward problem for elastodynamics formulated as a stochastic partial differential equation (SPDE) with random material properties and generalised forcing, accounting for applied loads and possible model discrepancies. In Section~\ref{sec:forwardgoveq}, we introduce the governing equations and model both the material properties and forcing as Gaussian random fields. In Section~\ref{sec:forwardgenmat}, we construct these random fields using an SPDE-based approach, which induces the desired covariance structure. We then discretise the problem in space using finite elements in Section~\ref{sec:forwardfem} and in time using an explicit scheme in Section~\ref{sec:forwardtime}. In Section~\ref{sec:forwardmodel}, we derive approximate expressions for the mean and covariance of the solution using first-order perturbation. Finally, in Section~\ref{sec:oscillator}, we demonstrate the accuracy of the proposed formulation through a single-degree-of-freedom example.
%--------------------------------------------------------------------------------          
\subsection{Governing equations \label{sec:forwardgoveq}}
%--------------------------------------------------------------------------------

The  elastodynamic equations for a solid body occupying a domain $\Omega \subset \mathbb{R}^d$, $d\in\{1,2,3\}$, over the time interval $t \in [0, T]$  are given by 
\begin{subequations}
    \begin{align}
    \begin{alignedat}{2}
    \vec{\nabla} \cdot \vec{\sigma}(\vec{u}) + \vec f &=\rho\ddot{\vec{u}}  &\qquad  &\text{in } \Omega \, \times \,  (0,T],\\
    \vec{u}&=\vec{0} & \qquad &\text{on }\Gamma_D \,  \times \, (0,T],\\
    \vec{\sigma}(\vec{u})\cdot \vec{n}&= \vec 0 &\qquad  &\text{on }\Gamma_N \, \times \, ( 0,T],\\
    \vec{u}(\vec{x},0)=\vec{u}_0\text{, } \dot{\vec{u}}(\vec{x},0)&=\dot{\vec{u}}_0  &\qquad &\text{for } \vec{x}\in\Omega .
    \end{alignedat}
\end{align}
\label{eq:goveq}
\end{subequations}
Here, $\vec{\sigma}$ is the stress tensor, $\vec{u}$ the displacement vector, ${\vec{f}}$ the generalised body force, $\rho$ the density, $\ddot{\vec u}$ the acceleration, $\vec u_0 $ the initial displacement, and $\dot{\vec u}_0 $ the initial velocity.  Without loss of generality, we consider only homogeneous boundary conditions on the Dirichlet boundary~$\Gamma_D$ and the Neumann boundary~$\Gamma_N$ with outward normal~$\vec{n}$.
We assume the linear elastic constitutive relation   
\begin{equation}
    \vec{\sigma}(\vec{u}) = \mathbb{C} :\vec{\epsilon} (\vec{u}),
\end{equation}
where $\mathbb{C}$ is the fourth-order tensor and~$\vec{\epsilon}$ is the strain tensor,   
\begin{equation}
    \vec{\epsilon} (\vec{u}) = \frac{1}{2}\left( \nabla \vec{u} + \left ( \nabla \vec{u} \right )^\trans \right).
\end{equation}

The generalised body force~$\vec  f$ and the constitutive tensor~$\mathbb C$ are random.  The constitutive tensor depends on a scalar random material property field $\kappa$, i.e. \mbox{$\mathbb{C}(\vec{x})=\mathbb{C}(\kappa(\vec{x}))$}. Here, $\kappa$ is a spatial Gaussian random field
\begin{equation}
    \kappa(\vec{x}) \sim \mathcal{GP}\left(\overline{\kappa}(\vec{x}),{c}_{\kappa}(\vec{x},\vec{x}^\prime)\right)\, ,    
\end{equation}
with mean $\overline{\kappa}$ and spatial covariance ${c}_{\kappa}$, 
\begin{subequations}\label{eq:kappa_moments}
\begin{align}
	\overline{\kappa}(\vec{x}) &= \expect{\left[\kappa(\vec{x})\right]} \, ,\\
       	{c}_{\kappa}(\vec{x},\vec{x}^\prime) & = \cov{\left(\kappa(\vec{x}),\kappa(\vec{x}^{\prime})\right)} = 
    \expect{\left[\big(\kappa(\vec{x})-\overline{\kappa}(\vec{x})\big)\big(\kappa(\vec{x}^\prime)-\overline{\kappa}(\vec{x}^\prime)\big)\right]}.\, 
\end{align}
\end{subequations}
The forcing $\vec{f}$ is a spatio-temporal Gaussian random field
\begin{equation}\label{eq:force}
	\vec{f}(\vec{x},t) \sim \mathcal{GP}\left(\overline{\vec{f}}(\vec{x},t),\vec{c}_{\vec{f}}(\vec{x},\vec{x}^\prime)\delta(t-t^\prime)\right)  ,
\end{equation}
with mean $\overline{\vec{f}}$ and spatial covariance $\vec{c}_{\vec{f}}$, 
\begin{subequations}\label{eq:force_moments}
\begin{align}
	\overline{\vec{f}}(\vec{x},t) &= \expect{\left[\vec{f}(\vec{x},t)\right]} \, ,\\
    \vec{c}_{\vec{f}}(\vec{x},\vec{x}^\prime) & = \cov{\left(\vec{f}(\vec{x},t),\vec{f}(\vec{x}^{\prime},t)\right)} = 
    \expect{\left[\left(\vec{f}(\vec{x},t)-\overline{\vec{f}}(\vec{x},t)\right)\left(\vec{f}(\vec{x}^\prime,t )-\overline{\vec{f}}(\vec{x}^\prime,t )\right)^{\trans}\right]} .
   \end{align}
\end{subequations}
The spatio-temporal covariance of~$\vec f$ is assumed to be separable in space and time, with $\vec{c}_{\vec{f}}$ describing spatial and the Dirac delta function $\delta$ temporal correlations. This yields a spatially correlated, temporally white noise representation of the force. 
%These spatially correlated random fields are generated by solving linear operators driven by white noise, as detailed in the next subsection.

%       
%--------------------------------------------------------------------------------          
\subsection{SPDE formulation of random fields \label{sec:forwardgenmat}}
%--------------------------------------------------------------------------------
%
Spatial random fields can be constructed as solutions of SPDEs driven by Gaussian noise. Specifically, consider 
\begin{equation}\label{eq:genspde}
    \set{L} s(\vec{x}) = \xi(\vec{x}), 
\end{equation}
where $\set{L}$ is a differential operator and $\xi$ is the Gaussian white noise field
\begin{equation}
	\xi (x) \sim \set{GP} \left ( 0, \delta(\vec x - \vec {x}^\prime) \right) .
\end{equation}
The operator $\set{L}$ determines the correlation structure of the random field. In this work, we choose
\begin{equation}\label{eq:matspde}
	\set{L} \coloneqq  \tau \left( \eta^2 - {\vec \nabla}^2 \right)^\beta, 
\end{equation}
where ${\vec \nabla}^2$ is the Laplacian. The parameter $\tau$ controls the variance $\sigma$, the parameter $\eta$ represents the inverse of the correlation length $l$ and the parameter $\beta$  defines the smoothness of the resulting random field. The field $s$ represents either the random material property field $\kappa$ or one of the spatial components~$f_i$ of the forcing field. For simplicity, we do not consider cross-correlations between different forcing components; however, they could be incorporated by generalising the operator~$\set{L}$ to a matrix operator. 

The SPDE~\eqref{eq:genspde} defines the Gaussian random field 
\begin{equation}
	s(\vec x) \sim \set{GP} \left ( 0, c_s(\vec x, \vec{x}^\prime )\right )
\end{equation}
with covariance~$c_s$. The solution of SPDE can be written as the convolution of its Green's function
$G$ with the driving noise~$\xi$,
\begin{equation}
    s(\vec{x}) = \int_\Omega G(\vec{x},\vec{y})\,\xi(\vec{y})\,\D\vec{y}.
\end{equation}
The covariance of $s$ is therefore given by 
\begin{equation}
    c_s(\vec{x},\vec{x}^{\prime}) = \cov{(s(\vec{x}),s(\vec{x}^{\prime}))} = \expect{\left[s(\vec{x})s(\vec{x}^{\prime})\right]} = \int_{\Omega}\int_{\Omega} G(\vec{x},\vec{y})\,\cov{\left(\xi(\vec{y}),\xi(\vec{y}^{\prime})\right)}\,G(\vec{x}^{\prime},\vec{y}')\,\D\vec{y}\,\D\vec{y}'.
\end{equation}
Since $\cov{\left(\xi(\vec{y}),\xi(\vec{y}^{\prime})\right)} = \delta(\vec y - \vec {y}^\prime)$, this simplifies to  
\begin{equation}
    c_s(\vec{x},\vec{x}^{\prime}) = \int_\Omega G(\vec{x},\vec{y})\,G(\vec{x}^{\prime},\vec{y})\,\D\vec{y}. 
\end{equation}
Thus, the covariance~$c_s$ is induced directly by the operator $\set{L}$. 

For completeness, the operator given in~\eqref{eq:matspde} leads to Mat\'ern random fields; the details of the SPDE–Mat\'ern correspondence are provided in \ref{sec:matern}. 

%--------------------------------------------------------------------------------
\subsection{Spatial finite element discretisation \label{sec:forwardfem}}
%--------------------------------------------------------------------------------
%
The elastodynamic problem~\eqref{eq:goveq} and the SPDE~\eqref{eq:genspde} are discretised using a standard finite element approach after establishing their weak forms. The problem domain~$\Omega$ is partitioned into a mesh with~$n_e$ elements and~$n_n$ nodes. To distinguish continuum fields from their discrete counterparts, continuum vectors are denoted by bold italic symbols, whereas algebraic arrays and matrices by upright bold symbols. 

The displacement field is interpolated as
\begin{equation}
    \vec{u}(\vec{x},t) \approx \sum_{{\mathrm{j}}=1}^{n_u} \phi_{\mathrm{j}}(\vec{x})\,\ary{u}_{\mathrm{j}}(t).
\end{equation}
where $\phi_{\mathrm{j}}$ are the basis functions associated with the~$n_u$ non-Dirichlet nodes and $\ary{u}_{\mathrm{j}}$ the corresponding nodal values. The semi-discrete elastodynamic equation with viscous damping is
\begin{equation}\label{eq:semi-discrete}
    \mat{M}\,\ddot{\ary{u}}
    + \mat{D}\,\dot{\ary{u}}
    + \mat{K}(\ary{\kappa})\,\ary{u}
    = \ary{f},
\end{equation}
where $\mat{M}$ is the mass matrix, $\mat{K}$ the stiffness matrix, $\mat{D}$ the damping matrix, and $\ary{f}$ the nodal force vector. The damping matrix is typically taken as a linear combination of $\mat{M}$ and $\mat{K}$ (e.g. Rayleigh damping). 

The same element mesh is used to discretise the SPDE~\eqref{eq:genspde} for the random material property and the forcing fields. The material property is approximated as 
\begin{equation}
    \kappa(\vec{x}) \approx \sum_{{\mathrm{j}}=1}^{n_n} \phi_{\mathrm{j}}(\vec{x})\, {\kappa}_{\mathrm{j}}.
\end{equation}
The summation is over all mesh nodes because homogeneous Neumann boundary conditions are assumed.  The discretisation yields the system of equations
\begin{equation}\label{q:RFSPDE-kappa}
\mat{L}_{\kappa}\,\ary{\kappa} = \mat{M}_{\kappa}^{1/2}\ary{\xi}, \qquad \ary{\xi}\sim\mathcal{N}(\ary{0},\mat{I}), 
\end{equation}
where~$\mat{M}_{\kappa}$ is the diagonal lumped Gram matrix, $\mat{L}_\kappa$ is the finite element discretised operator from~\eqref{eq:matspde} (cf.~\eqref{eq:discrete-lap}) and~$\ary \xi$ is an array of independent standard Gaussian variables. By the linear transformation property of Gaussian densities,~$\ary \kappa$ has the covariance 
\begin{equation}\label{q:RFSPDE-kappa-cov}
	\mat{C}_{\kappa} = \mat{L}_{\kappa}^{-1}\,\mat{M}_{\kappa}\,\mat{L}_{\kappa}^{-\trans}.
\end{equation}
Hence, the material property array has the multivariate Gaussian density 
\begin{equation}
	\ary \kappa \sim \set N \left (\overline{\ary \kappa} , \ary C_{\kappa} \right ) .
\end{equation}
The mean~$\overline{\ary \kappa}$  is obtained by evaluating~$\overline{\kappa} (\vec x)$ at the mesh nodes. In the subsequent formulation, the element-wise material property is obtained by averaging the nodal values over each element. The corresponding covariance is obtained by projecting the nodal covariance~$\ary C_{\kappa}$ onto the element space. With slight abuse of notation, we use the same symbol~$\ary C_{\kappa}$ for both the nodal and element covariance.

The spatial covariance of the force array~$\ary f$  is obtained by discretising the respective SPDE using the same approach as before. As introduced in Section~\ref{sec:forwardgoveq}, the random forcing~$\vec f$ is spatially correlated and temporally white. The discretised SPDE reads 
\begin{subequations}\label{q:RFSPDE-force}
\begin{align}
 &\mat{L}_f\,(\mat{M}_f^{-1}\ary{f}) = \mat{M}_f^{1/2}\ary{\xi},
    \qquad
    \ary{\xi} \sim \mathcal{N}(\ary{0},\mat{I}),
\intertext{yielding the spatial covariance}
&\mat{C}_f = \mat{M}_f\,\mat{L}_f^{-1}\,\mat{M}_f\,\mat{L}_f^{-\trans}\,\mat{M}_f^{\trans}.
\end{align}
\end{subequations}
The two additional Gram matrices~$\ary {M}_f$, in comparison to~\eqref{q:RFSPDE-kappa-cov}, arise because~$\ary f$ are the nodal forces obtained by integrating the field~$\vec f$. Hence, the force array has the multivariate Gaussian density 
\begin{equation}
	\ary f \sim \set N \left (\overline{\ary f} , \ary C_{f}  \delta(t - t^\prime) \right ).
\end{equation}
%
%Equation~\eqref{eq:semi-discrete} represents the spatially discretized stochastic elastodynamic problem incorporating random stiffness and random loads. The subsequent section describes the temporal discretization and the stochastic time integration scheme used to propagate these uncertainties.

%
%--------------------------------------------------------------------------------          
\subsection{Temporal discretisation \label{sec:forwardtime}}
%--------------------------------------------------------------------------------
%
We consider the temporal discretisation of the semi-discrete elastodynamic equation~\eqref{eq:semi-discrete} with a random material property and forcing fields. The semi-discrete equation expressed in state-space form reads
\begin{equation}\label{eq:EOM_SP_Syst}
	\D\,
	\begin{pmatrix}
		{\ary{u}}\\
		\dot{\ary{u}}
	\end{pmatrix}=
	\begin{pmatrix}
		\mat{0} & \mat{I}\\
		- \mat{M}^{-1}\mat{K}(\ary{\kappa}) & -\mat{M}^{-1}\mat{D} 
	\end{pmatrix}
	\begin{pmatrix}
		\ary{u}\\
		\dot{\ary{u}}
	\end{pmatrix}
   \D t
	+
	\begin{pmatrix}
		\ary{0}\\
		{\mat{M}^{-1}}
	\end{pmatrix}
	\ary{f}\D t \, .
\end{equation}
Defining the state array 
\begin{equation}
\ary{v} = \begin{pmatrix} \ary{u}^{\trans} & \dot{\ary{u}}^{\trans} \end{pmatrix}^{\trans},
\end{equation}
the equation can be written as the (It\^{o}) stochastic differential equation 
\begin{equation}\label{eq:EOM_SP}
		\D{\ary{v}}    
		=\mat{F}(\ary{\kappa})\ary{v} \D t+\mat{G}\overline{\ary{f}}\D t+\mat{G}\D \ary{\beta} \, ,
\end{equation}
where $\mat{F}(\ary{\kappa})$ is the state matrix, $\mat{G}$ the input matrix, and $\ary{\beta}$ a spatially correlated Brownian motion. Its increments are Gaussian with zero mean and covariance proportional to the time increment, i.e.,
\begin{equation} \label{eq:definition_Brownian}
	\ary \beta(t) - \ary \beta(t')\sim \set N(\ary 0, \ary C_f (t-t^\prime)), \quad  0 \le t^\prime < t . 
\end{equation}

A stochastic time-marching scheme is required for the discrete solution of~\eqref{eq:EOM_SP} as it involves an It\^{o} integral for the Brownian forcing. We use the stochastic Verlet algorithm in its position-based leapfrog form, which is quasi-symplectic and second-order accurate for elastic forces. In the absence of damping, the state evolution is identical to the classical explicit central-difference scheme or the explicit Newmark scheme, except for the treatment of the spatially correlated Brownian motion contribution.

After partitioning the time~$(0, T]$ into equidistant intervals of size~$\Delta t = t_{n+1} - t_n$, with~$0< t_1< t_2 < \dotsc \le T$, the Verlet state evolution (see~\ref{sec:app-verlet}) can be expressed as 
\begin{subequations}\label{eq:vvFEM}
\begin{align}
 \ary{v}_{n+1} &= \mat{A}(\ary{\kappa})\,\ary{v}_n 
    + \Delta t\,\mat{B}\,\overline{\ary{f}}_n + \ary{\zeta}_n ,\\[3pt]
\intertext{with}
\mat{A}(\ary{\kappa}) &=
    \begin{pmatrix}
        \mat{I}-\tfrac{\Delta t^2}{2}\mat{M}^{-1}\mat{K}(\ary{\kappa}) & \Delta t \left( \mat{I}-\tfrac{\Delta t}{2} \mat{M}^{-1}\mat{D}-\tfrac{\Delta t^2}{4}\mat{M}^{-1}\mat{K}(\ary{\kappa})\right)\\
        -\Delta t \mat{M}^{-1}\mat{K}(\ary{\kappa}) & \mat{I}-\Delta t \mat{M}^{-1}\mat{D}-\tfrac{\Delta t^2}{2}\mat{M}^{-1}\mat{K}(\ary{\kappa})
    \end{pmatrix},
    \quad
    \mat{B} =
    \begin{pmatrix}
        \tfrac{\Delta t }{2}\mat{M}^{-1} \\
        \mat{M}^{-1}
    \end{pmatrix}, 
    \quad 
    \ary{\zeta}_n =\mat{B}\,\Delta\ary{\beta}_n  \, . 
\end{align}
\end{subequations}
The Brownian component~$\ary{\zeta}_n$ has the multivariate Gaussian density 
\begin{subequations}
\begin{align}
	 \ary{\zeta}_n &\sim \set N (\ary 0, \ary C_{\zeta}),  \\
	 \intertext{which has according to~\eqref{eq:definition_Brownian} the covariance}
	  \mat{C}_{\zeta} & = \mat{B}\,(\mat{C}_f\,\Delta t)\,\mat{B}^\trans \,. 
\end{align}
\end{subequations}
%

% In the introduced scheme, stability is guaranteed in the absence of damping,  provided the time step satisfies the usual central-difference condition
% %
% \begin{equation}
% 	\Delta t \leq 2/\omega_{\max} , 
% \end{equation}
% %
% where $\omega_{\max}$ is the maximum frequency of the discretised system. 
The state evolution~\eqref{eq:vvFEM} provides an explicit relation for propagating the state vector under material and forcing uncertainties. It depends linearly on the forcing arrays~$\overline{\ary{f}}$  and~$\ary{\zeta}_n$, and nonlinearly on the material property array~$\ary \kappa$.  Although the inputs are Gaussian, the state evolution is nonlinear in~$\ary \kappa$, so that the probability density of the state is generally non-Gaussian and has no closed form.  

%--------------------------------------------------------------------------------
\subsection{Perturbation-based stochastic solution  \label{sec:forwardmodel}}
%--------------------------------------------------------------------------------
%
The Verlet state evolution~\eqref{eq:vvFEM} defines a mapping~$\ary{v}_{n+1}(\ary{v}_{n},\ary \kappa, \ary{\zeta}_n)$ with the previous state~$\ary{v}_{n}$, material property~$\ary \kappa$, and the Brownian increment~$\ary{\zeta}_n$ as its arguments.  A first-order expansion of this mapping about the mean values~$\overline{\ary{v}}_{n}$, $\overline {\ary \kappa}$ and~$\overline{\ary{\zeta}}_n$ yields 
\begin{equation}
\begin{aligned}
	\ary v_{n+1} (\ary{v}_{n},\ary \kappa, \ary{\zeta}_n)  & \approx \ary v_{n+1} (\overline{\ary{v}}_{n}, \overline{\ary \kappa}, \overline{\ary{\zeta}}_n)
  + \frac{\partial \ary{v}_{n+1}}{\partial \ary{v}_{n}}\Big|_{(\overline{\ary{v}}_{n}, \overline{\ary \kappa}, \overline{\ary{\zeta}}_n)} (\ary{v}_{n}-\overline{\ary{v}}_n) 
 + \frac{\partial \ary{v}_{n+1}}{\partial \ary{\kappa}}\Big|_{(\overline{\ary{v}}_{n}, \overline{\ary \kappa}, \overline{\ary{\zeta}}_n)    }(\ary{\kappa}-\overline{\ary{\kappa}})  + \frac{\partial \ary{v}_{n+1}}{\partial \ary{\zeta}_n}\Big|_{(\overline{\ary{v}}_{n}, \overline{\ary \kappa}, \overline{\ary{\zeta}}_n)}\,\ary{\zeta}_n  \, . 
\end{aligned}
\end{equation}
The constant term evaluates to
\begin{equation}
	\overline{\ary v}_{n+1} = \ary v_{n+1} (\overline{\ary{v}}_{n}, \overline{\ary \kappa}, \overline{\ary{\zeta}}_n)  = \overline{\mat{A}} \, \overline{\ary{v}}_n  + \Delta t\,\mat{B}\,\overline{\ary{f}}_n \, ,
    \label{eq:mean-update}
\end{equation}
where~$\overline{\mat A } = \mat{A}(\overline{\ary{\kappa}})$ and~$\overline{\ary{\zeta}}_n = \ary 0$. The true mean $\overline{\ary v}_{n+1}$ involves the mean~$\overline{\mat{A}(\ary{\kappa})\, \ary{v}_n}$. Taking instead~$\mat{A}(\overline{\ary{\kappa}})\, \overline{\ary{v}}_n$ constitutes a first-order approximation. After evaluating the gradients the linearised state evolution can be written as 
\begin{subequations}\label{eq:fluc-update}
\begin{align}
\ary v_{n+1} &=  \overline{\ary v}_{n+1}   + \overline{\mat A} (\ary{v}_{n}-\overline{\ary{v}}_n) + \mat{J}_n (\ary{\kappa}-\overline{\ary{\kappa}})  + \ary{\zeta}_n  \, ,
\\
\intertext{with}
\mat{J}_n  (\ary{\kappa}-\overline{\ary{\kappa}}) 
&:= \sum_{\alpha=1}^{n_{el}} 
\left ( \frac{\partial \mat{A}}{\partial \kappa_\alpha}\Big|_{\overline{\ary{\kappa}}}\,
\overline{\ary{v}}_n \right ) \, (\ary{\kappa}-\overline{\ary{\kappa}})  .
\end{align}
\end{subequations}

We aim to compute the joint probability density~$p(\ary v_{n+1}, \ary \kappa)$ required for Bayesian filtering.  All input random variables are Gaussian, i.e., 
\begin{equation}
	\ary \kappa  \sim \set N(\overline{\ary \kappa}, \ary{C}_{\kappa}), \quad  \ary{\zeta}_n \sim \set N (\ary 0, \ary C_{\zeta}) \, .
\end{equation}
Since the linearised state evolution~\eqref{eq:fluc-update} represents a linear transformation of Gaussian variables, the evolved state is Gaussian as well. Using standard identities for covariance,  we can derive
\begin{subequations}
\begin{align}
\cov(\ary{v}_{n+1},\ary{v}_{n+1}) &= \overline{\mat{A}}\,\cov(\ary{v}_{n},\ary{v}_{n})  \,\overline{\mat{A}}^\trans
+ \overline{\mat{A}}\, \cov(\ary{v}_{n},\ary{\kappa})  \,\mat{J}_n^\trans
+  \mat{J}_n  \, \cov(\ary{\kappa},\ary{v}_{n})  \,   \overline{\mat{A}}^\trans 
+\mat{J}_n\, \cov(\ary \kappa, \ary \kappa) \,\mat{J}_n^\trans
+\cov( \ary {\zeta}_n , \ary{\zeta}_n) , \\[4pt]
\cov(\ary{v}_{n+1},\ary{\kappa}) &= \overline{\mat{A}}\,  \cov(\ary{v}_{n},\ary{\kappa})
+\mat{J}_n\, \cov(\ary{\kappa},\ary{\kappa}) . 
\end{align}
\end{subequations}
Or, more compactly, 
\begin{subequations}\label{eq:cov-compact}
\begin{align}
\mat{C}_{n+1} &= \overline{\mat{A}}\,\mat{C}_{n}\,\overline{\mat{A}}^\trans
+  \overline{\mat{A}}\,\mat{C}_{n,\kappa}\,\mat{J}_n^\trans
 + \mat{J}_n \,\mat{C}_{\kappa, n}\,\overline{\mat{A}}^\trans
+\mat{J}_n\,\mat{C}_\kappa\,\mat{J}_n^\trans
+\mat{C}_{\zeta}, \label{eq:cov-state}\\[4pt]
\mat{C}_{n+1,\kappa} &= \overline{\mat{A}}\,\mat{C}_{n,\kappa}
+\mat{J}_n\,\mat{C}_\kappa . \label{eq:cov-cross}
\end{align}
\end{subequations}
Although the Brownian increments $\ary{\zeta}_n$ are spatially correlated, the forcing is temporally white. Consequently, $\ary{\zeta}_n$ is independent of all quantities determined up to time $t_n$, in particular of $\ary{v}_n$ and $\ary{\kappa}$.  Furthermore, the initial state~$\ary v_0$ is usually independent of the material field~$\ary \kappa$ so that~$\mat C_{0,\kappa} = \ary 0$. The second equation then gives~$\mat{C}_{1,\kappa} = \mat J_0  \mat C_{\kappa}$, after which the recursion applies.

Finally, the joint density of state and material property at time~$t_{n+1}$ is
\begin{equation}
p ( \ary{v}_{n+1}, \ary{\kappa})
=
\mathcal{N}\!\left(
\begin{pmatrix}
\overline{\ary{v}}_{n+1}\\
\overline{\ary{\kappa}}
\end{pmatrix},
\begin{pmatrix}
\mat{C}_{n+1} & \mat{C}_{n+1,\kappa}\\
\mat{C}^\trans_{n+1,\kappa} & \mat{C}_{\kappa}
\end{pmatrix}
\right).
\label{eq:probfm}
\end{equation}
As mentioned, the true density is generally not Gaussian. This approximation is valid only when the variance of the material property and the time-step size are small.

%--------------------------------------------------------------------------------
\subsection{Illustrative example \label{sec:oscillator}}
%--------------------------------------------------------------------------------
To verify the accuracy of the perturbation-based stochastic formulation, we assess the predictive capability of the perturbation-based stochastic solution using a simple single-degree-of-freedom (SDOF) system and benchmark the results against Monte Carlo simulations. We consider the SDOF system subject to Gaussian white noise forcing,
\begin{equation}
    \ddot{u}(t) = -\omega^2 u(t) - \frac{\gamma}{m} \dot{u}(t) +  \frac{1}{m}{(\overline{f}(t)+\xi(t))},
\end{equation}
where $\omega=\sqrt{k/m}$ is the natural frequency with $k$ being the stiffness and $m$ being the mass, $\gamma$ is the damping coefficient, $\overline{f}(t)$ is the deterministic component of the force, and $\xi(t)$ is zero-mean white noise with spectral density $\sigma^2$, which is the stochastic component of the force. 

For numerical validation, we choose natural frequency of $\omega_s = 10$, a damping coefficient $\gamma = 1.0$, and a spring constant $k \sim \mathcal{N}(k_0, \sigma_k^2)$ that remains constant in time. The mean and standard deviation of the spring constant are $k_0 = 100$ and $\sigma_k = 0.05k_0$, respectively. The random external forcing consists of two harmonic components, $f(t) = \overline{f}(t) + \xi(t)$, with $\overline{f} = f_0[\sin(0.31\omega_s t) + \sin(0.62\omega_s t)]$, and $\xi(t)\Delta t = \Delta \beta(t) \sim N(0,\sigma_f^2\Delta t)$ with $\sigma_f = 0.05f_0$ and $f_0=1.0$. The integration time step is chosen as $\Delta t / T_s = 3.2 \times 10^{-3}$, where $T_s = 2\pi / \omega_s$ denotes the natural period.

\begin{figure}[!t]
	\centering
	\subfloat[Displacement samples]{\includegraphics[width=0.49\textwidth]{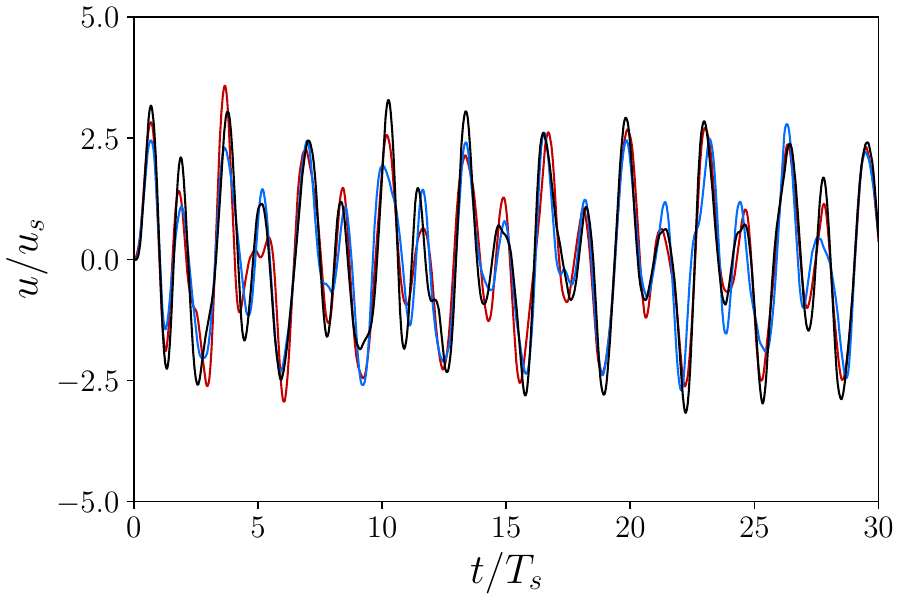}}\hfill	
    \subfloat[Standard deviation]{\includegraphics[width=0.49\textwidth]{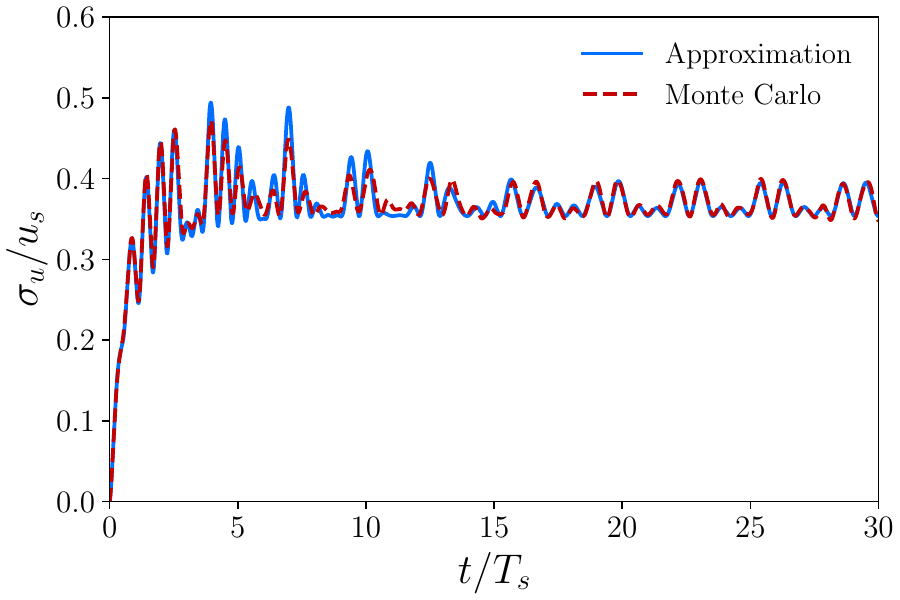}}
	\caption{
Accuracy of the perturbation-based stochastic model for a single-degree-of-freedom (SDOF) oscillator compared with Monte Carlo (MC) results. The figures show the displacement ($u$) and standard deviation of the solution ($\sigma_u$), normalised by the static displacement $u_s = f_0 / \overline{k}$.  
}
	\label{fig:SDOF_Forward}
\end{figure}

\begin{figure}[!t]
    \centering
    \subfloat[Varying material uncertainty \label{fig:SDOF_Convergence_STD}]
        {\includegraphics[width=0.48\textwidth]{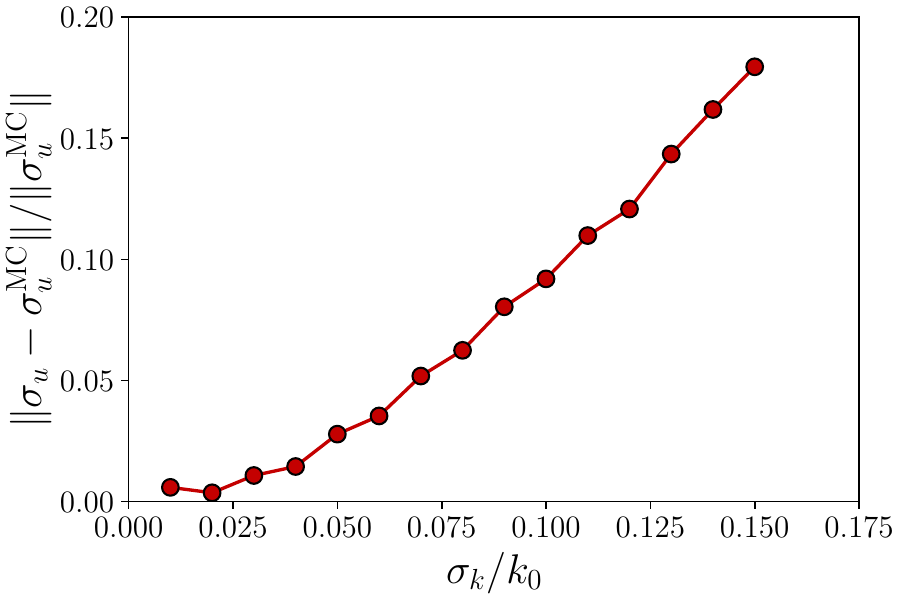}\vspace{-1mm}}\hfill
    \subfloat[Varying time step \label{fig:SDOF_Convergence_TS}]
        {\includegraphics[width=0.48\textwidth]{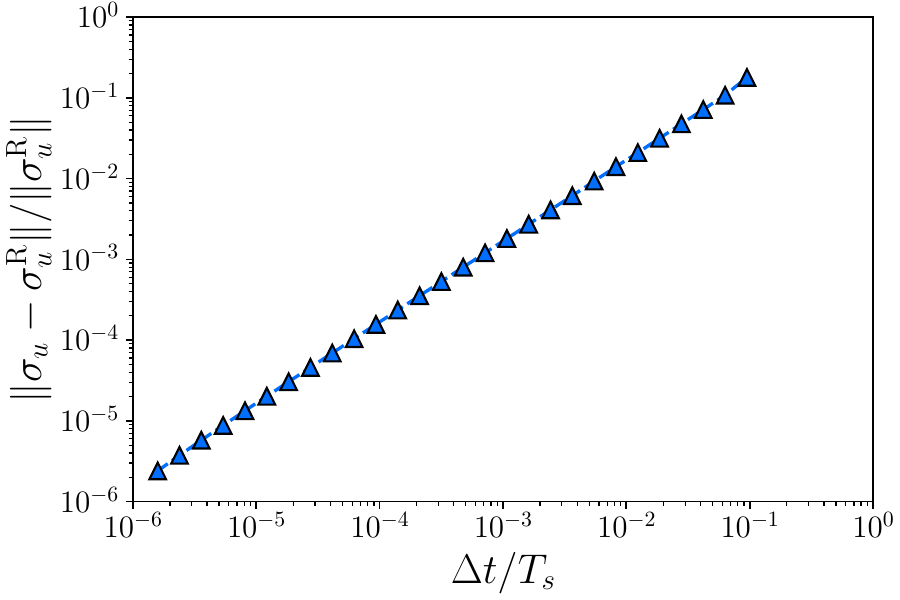}\vspace{-1mm}}
    \caption{Convergence plots for the relative perturbation approximate standard deviation $\sigma_u$ for the single-degree-of-freedom problem with material uncertainty. In (a), only the standard deviation of the material uncertainty $\sigma_k$ is varied. The empirical standard deviation $\sigma_u^{\mathrm{MC}}$ is obtained by Monte Carlo sampling and considered as exact. Similarly, in (b), the standard deviation errors are shown for different time step sizes, taking the smallest time step solution, $\sigma_u^R$, as the reference.}
    \label{fig:SDOF_Convergence}
\end{figure}

Figure~\ref{fig:SDOF_Forward} shows a few samples of displacement as a function of time and compares the predictive variance obtained from the forward model with Monte Carlo estimates, normalized by the static displacement $u_s = f_0 /k_0$. We notice that the predictive standard deviation of the forward model matches quite well with Monte Carlo estimates. To further evaluate the accuracy of the predictive variance, we examine different magnitudes of material uncertainty by varying $\sigma_k/k_0$, as shown in Fig.~\ref{fig:SDOF_Convergence_STD}. Good agreement with the Monte Carlo estimates is observed for small variances, while the error increases with larger $\sigma_k/k_0$ owing to the first-order accuracy of the expansion. Figure~\ref{fig:SDOF_Convergence_TS} additionally demonstrates the convergence of the variance with decreasing time step.

\begin{comment}
\begin{figure}[!t]
	\centering
	\subfloat[Displacement]{\includegraphics[width=0.49\textwidth]{figs/forward/onedof_mu_mean_with_samples_both.pdf}}\hfill	
    \subfloat[Standard deviation]{\includegraphics[width=0.49\textwidth]{figs/forward/onedof_mu_std_prediction_vs_mc_both.pdf}}
	\caption{
Accuracy of the perturbation-based stochastic model for a single-degree-of-freedom (SDOF) oscillator compared with Monte Carlo (MC) results. The figures show the displacement ($u$) and standard deviation of the solution ($\sigma_u$), normalised by the static displacement $u_s = f_0 / \overline{k}$.  
}
	\label{fig:SDOF_Forward}
\end{figure}
\end{comment}

%% file: filtering.tex
%
%--------------------------------------------------------------------------------          
\section{Statistical Finite Elements through Bayesian Filtering \label{sec:bayesfil}}
%--------------------------------------------------------------------------------
%

In this section, we employ Bayesian inference to incorporate observational data and update the solution of the elastodynamic forward problem. In Section \ref{sec:intbayesfil}, we introduce the statistical observation model and formulate the Bayesian inference as a sequential filtering procedure consisting of prediction and update steps. Subsequently, in Sections \ref{sec:bayespred}, we derive expressions for the predicted and updated probability densities of the state vector. As noted in Section \ref{sec:forwardmodel}, we augment the state vector to include the material field in addition to the positions and velocities of the forward problem. Finally, in Secton \ref{sec:bayes-params}, we address the estimation of the covariance hyperparameters.

%--------------------------------------------------------------------------------
\subsection{Bayesian filtering framework \label{sec:intbayesfil}}
\begin{figure}[!t]
    \centering  
       {\includegraphics[width=0.75\textwidth]{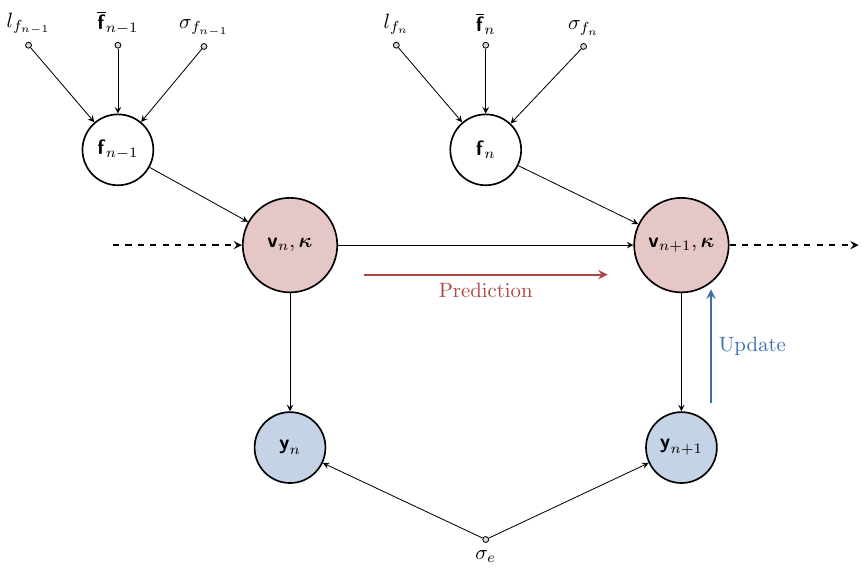}\vspace{-1mm}}
    \caption{Graphical model of statFEM for solid dynamics with augmented-state filtering. The large circles represent random variables that are observed or inferred. The augmented state vector $\left(\ary{v}_{n}^\trans\; \ary{\kappa}^\trans\right)^\trans$ is updated sequentially as observations $\ary{y}_n$ are assimilated. The hyperparameters $\sigma_f$ and $l_f$ correspond to the standard deviation and correlation length of the model misspecification field. The measurement noise is independent and identically distributed with standard deviation $\sigma_e$.}
    \label{fig:bayesian_framework}
\end{figure}

%--------------------------------------------------------------------------------
We are given the observations
\begin{equation}
\ary{y}_{1:N}=\{\ary{y}_1,\dots,\ary{y}_{N}\},
\end{equation}
along with the probabilistic forward model~\eqref{eq:probfm}. The observation model is assumed to be
\begin{align}\label{eq:linear-sys}
\ary{y}_{n+1} &= \mat{H}\,\ary{v}_{n+1} + \ary{e}_{n+1}, 
\qquad \ary{e}_{n+1} \sim \mathcal{N}(\ary{0},\,\mat{C}_e),
\end{align}
where $\mat{H}$ is the observation matrix and $\mat{C}_e$ denotes the measurement noise covariance. The aim is to infer the posterior distribution $p(\ary{v}_{n+1},\ary{\kappa}|\ary{y}_{1:n+1})$ of the current augmented state given all observations up to time $t_{n+1}$. The assumed Markovian statistical model, which encodes the conditional independence structure of the random variables, is shown in Fig.~\ref{fig:bayesian_framework}. This sequential Markovian structure allows the posterior to be computed recursively from the previous time step using two stages: prediction and update.

The prediction step propagates the uncertainty in the state from $t_n$ to $t_{n+1}$ through the transition model described in Section~\ref{sec:forwardmodel}. Starting from the initial probability density $p(\ary{v}_{0},\ary{\kappa})$, the prediction step determines the prior
\begin{equation}\label{eq:predict-bayes-general}
p(\ary{v}_{n+1},\ary{\kappa}|\ary{y}_{1:n})
= \int p(\ary{v}_{n+1},\ary{\kappa},\ary{v}_n,{\ary{f}_n}|\ary{y}_{1:n})\,\D\ary{v}_n \, {\D\ary{f}_n}= \int p(\ary{v}_{n+1}|\ary{v}_n,{\ary{\kappa}},{\ary{f}_n},\ary{y}_{1:n})\,p(\ary{v}_n,{\ary{\kappa}}|{\ary{f}_n},\ary{y}_{1:n}) {p(\ary{f}_n|\ary{y}_{1:n})}\,\D\ary{v}_n\,{\D\ary{f}_n}.
\end{equation}
According to the graphical model in Fig.~\ref{fig:bayesian_framework}, $\ary{v}_{n+1}$ is conditionally independent of the previous observations given $(\ary{v}_n,\ary{\kappa},\ary{f}_n)$. Moreover, since the forcing is temporally white, $\ary{f}_n$ is independent of $(\ary{v}_n,\ary{\kappa},\ary{y}_{1:n})$. Hence,
\begin{equation}\label{eq:predict-bayes}
p(\ary{v}_{n+1},\ary{\kappa}|\ary{y}_{1:n})
= \int p(\ary{v}_{n+1}|\ary{v}_n,{\ary{\kappa}},{\ary{f}_n})\,p(\ary{v}_n,{\ary{\kappa}}|\ary{y}_{1:n}) {p(\ary{f}_n)}\,\D\ary{v}_n\,{\D\ary{f}_n}.
\end{equation} 

The update step then updates prediction by applying Bayes’ rule to include the new observation:
\begin{equation}\label{eq:update-bayes-general}
p(\ary{v}_{n+1},\ary{\kappa}|\ary{y}_{1:n+1})
= \frac{p(\ary{y}_{n+1}|\ary{v}_{n+1},\ary{\kappa},\ary{y}_{1:n})\,p(\ary{v}_{n+1},\ary{\kappa}|\ary{y}_{1:n})}
       {p(\ary{y}_{n+1}|\ary{y}_{1:n})}.
\end{equation}
Since the observation depends on the {augmented} state only through $\ary{v}_{n+1}$ {and $\ary{\kappa}$} (see Fig.~\ref{fig:bayesian_framework}), the likelihood reduces to 
\begin{equation}\label{eq:update-bayes}
p(\ary{v}_{n+1},\ary{\kappa}|\ary{y}_{1:n+1})
= \frac{p(\ary{y}_{n+1}|\ary{v}_{n+1},\ary{\kappa})\,p(\ary{v}_{n+1},\ary{\kappa}|\ary{y}_{1:n})}
       {p(\ary{y}_{n+1}|\ary{y}_{1:n})},
\end{equation}
where the marginal likelihood in the denominator,
\begin{equation}\label{eq:marg-like-bayes}
p(\ary{y}_{n+1}|\ary{y}_{1:n})
= \int p(\ary{y}_{n+1}|\ary{v}_{n+1},\ary{\kappa})\,p(\ary{v}_{n+1},\ary{\kappa}|\ary{y}_{1:n})\,\D\ary{v}_{n+1}\,\D\ary{\kappa},
\end{equation}
ensures proper normalization. This completes the filtering equations at step $n+1$, which are used recursively in subsequent filtering steps.

%--------------------------------------------------------------------------------
\subsection{Predicted and updated probability densities \label{sec:bayespred}}
%--------------------------------------------------------------------------------
The derived filtering equations yield a natural state-augmented form, wherein the material property field  $\ary{\kappa}$ and state $\ary{v}_{n+1}$ are jointly inferred. The material property field is assumed to be without any dynamics, i.e. $\ary{\kappa}=\ary{\kappa}_{n+1}=\ary{\kappa}_{n}$ and is therefore updated only during the update step. Under the adopted linearisation, the joint predictive distribution in~\eqref{eq:predict-bayes} of the state and the material property field remains Gaussian. Specifically, the augmented vector
\(
\left(\ary{v}_{n+1}^\trans\; \ary{\kappa}^\trans\right)^\trans
\)
conditioned on the observations up to time $t_n$ has the joint probability density
\begin{equation}
p\,\left(
\begin{pmatrix}
\ary{v}_{n+1}\\
\ary{\kappa}
\end{pmatrix}
\Big|\,
\ary{y}_{1:n}
\right)
=
\mathcal{N}\,\left(
\begin{pmatrix}
\overline{\ary{v}}_{n+1}^{\, -}\\
\overline{\ary{\kappa}}_{n+1}^{\, -}
\end{pmatrix},
\begin{pmatrix}
\mat{C}_{n+1}^{\, -} & \mat{C}_{n+1,\kappa|n}^{\, -}\\
(\mat{C}_{n+1,\kappa|n}^{\, -})^\trans & \mat{C}_{\kappa|n}^{\, -}
\end{pmatrix}
\right).
\label{eq:joint-predictive}
\end{equation}
Here the superscript ${}^{-}$ denotes predictive (prior) quantities conditioned on the observations up to time $t_n$. Thus, $\overline{\ary{v}}_{n+1}^{\, -}$ and $\mat{C}_{n+1}^{\, -}$ are the predicted mean and covariance of the state at $t_{n+1}$ given the data $\ary{y}_{1:n}$ i.e. $p(\ary{v}_{n+1}|\ary{y}_{1:n})=\mathcal{N}(\overline{\ary{v}}_{n+1}^{\, -}, \mat{C}_{n+1}^{\, -})$, and the remaining terms follow analogously for the material property fields and cross–covariances.

The moments of the predictive distribution $p(\ary{v}_{n+1},\ary{\kappa}|\ary{y}_{1:n})$ in~\eqref{eq:predict-bayes} follow directly from the state update~\eqref{eq:fluc-update}, as detailed in Section~\ref{sec:predmean}. The resulting closed-form expressions, written in terms of the updated posterior statistics of the augmented state vector at time $t_n$, are summarized in Algorithm~\ref{alg:statFEMDyn}. The predictive mean and covariance of the state retain the same functional form as those of the forward model, cf.~\eqref{eq:mean-update} and~\eqref{eq:cov-compact}, with the distinction that the linearized dynamics matrices are evaluated at the posterior means, namely $\overline{\mat{A}}_n^{+}=\mat{A}(\overline{\ary{\kappa}}_n^{+})$ and $\mat{J}_n^{+}=\mat{J}_n(\overline{\ary{\kappa}}_n^{+},\overline{\ary{v}}_n^{+})$.

The joint posterior distribution of the state and material property fields also remain Gaussian. Specifically, the updated probability density of the augmented vector
\(
(\ary{v}_{n+1}^\trans\; \ary{\kappa}^\trans)^\trans
\)
conditioned on the observations up to time $t_{n+1}$ is given by
\begin{equation}
p\,\left(
\begin{pmatrix}
\ary{v}_{n+1}\\
\ary{\kappa}
\end{pmatrix}
\Big|\,
\ary{y}_{1:n+1}
\right)
=
\mathcal{N}\,\left(
\begin{pmatrix}
\overline{\ary{v}}_{n+1}^{\, +}\\
\overline{\ary{\kappa}}_{n+1}^{\, +}
\end{pmatrix},
\begin{pmatrix}
\mat{C}_{n+1}^{\, +} & \mat{C}_{n+1,\kappa|n+1}^{\, +}\\
(\mat{C}_{n+1,\kappa|n+1}^{\, +})^\trans & \mat{C}_{\kappa|n+1}^{\, +}
\end{pmatrix}
\right).
\end{equation}
It is obtained by applying Bayes’ rule~\eqref{eq:update-bayes} to combine the predictive distribution $p(\ary{v}_{n+1},\ary{\kappa}| \ary{y}_{1:n})$ with the likelihood $p(\ary{y}_{n+1}|\ary{v}_{n+1})$ (see~\ref{sec:app-kalman-update}). The complete prediction–update equations associated with the forward model developed in Section~\ref{sec:forwardmodel} are presented in Algorithm~\ref{alg:statFEMDyn} and are given in~\eqref{eq:kalman-pred} and~\eqref{eq:kalman-update}. Their recursive application propagates the joint mean and covariance of the augmented state, thereby providing a complete Bayesian update of the system state as new observations are assimilated.

In some applications, one may wish to treat $\ary{\kappa}$ as uncertain but not learnable from data or inferred through empirical Bayes (see Sec.~\ref{sec:bayes-params}). In such cases, the augmented state estimation can be simplified by suppressing the filtering update for $\ary{\kappa}$. Specifically, $\overline{\ary{\kappa}}_{n}^{+}$ and $\mat{C}_{\kappa|n}^{+}$ are held fixed at their prior values, and only the state statistics $(\overline{\ary{v}}_n^{+},\mat{C}_n^{+})$ are propagated and updated. The prediction and update steps for the state retain the same structure as in Algorithm~\ref{alg:statFEMDyn}, with the cross-covariance terms involving $\ary{\kappa}$ either omitted or treated as known.

\begin{algorithm}[H]
\caption{Bayesian prediction and update scheme}
\label{alg:statFEMDyn}
\begin{algorithmic}[1]
\Require
Initial state $\ary{v}_0$, prior statistics for material property field $\overline{\ary{\kappa}}_0^{\, +}$ and $\mat{C}_{\kappa|0}^{\, +}$, observations $\{\ary{y}_{n}\}_{n=1}^{N}$ and elastodynamic system matrices $\mat{A}(\cdot)$, $\mat{J}_n(\cdot,\cdot)$, $\mat{B}$, $\mat{H}$, $\mat{C}_\zeta$ and $\mat{C}_e$.
\vspace{0.5em}
\Ensure
Posterior statistics $\{\overline{\ary{v}}_n^{\, +},\overline{\ary{\kappa}}_n^{\, +}, \mat{C}_n^{\, +},\mat{C}_{\kappa|n}^{\, +},\mat{C}_{n,\kappa|n}^{\, +}\}_{n=0}^{N}$
\vspace{0.5em}

\Statex \textbf{Initialisation:} $\overline{\ary{v}}_0^{\, +} = \ary{v}_0$, $\mat{C}_0^{\, +} = \mat{0}$, $\mat{C}_{0,\kappa|0}^{\, +} = \mat{0}$, $\overline{\ary{\kappa}}_0^{\, +} = \ary{0}$ and $\mat{C}_{\kappa|0}^{\, +} = \mat{C}_{\kappa}$. 
\vspace{0.5em}
\Statex \textbf{Recursion:} Updated statistics at $t_n$ are known
\Statex  \hspace{1.2em} \textbf{Prediction step:}
{
\setlength{\jot}{6pt}   
\begin{subequations}\label{eq:kalman-pred}
\begin{align}
\overline{\ary{v}}_{n+1}^{\, -} 
&= \overline{\mat{A}}_n^{\,+}\,\overline{\ary{v}}_n^{\, +} 
+ \Delta t\,\mat{B}\,\overline{\ary{f}}_n,\\
\mat{C}_{n+1}^{\, -} 
&= \overline{\mat{A}}_n^{\,+}\,\mat{C}_n^{\, +}\,(\overline{\mat{A}}_n^{\,+})^\trans 
\;+\; \overline{\mat{A}}_n^{\,+}\,\mat{C}_{n,\kappa|n}^{\, +}\,(\mat{J}_n^{\, +})^\trans
\;+\; \mat{J}_n^{\, +}\,(\mat{C}_{n,\kappa|n}^{\, +})^{\trans}\,(\overline{\mat{A}}_n^{\,+})^\trans
\;+\; \mat{J}_n^{\, +}\,\mat{C}_{\kappa|n}^{\,+}\,(\mat{J}_n^{\, +})^\trans
\;+\; \mat{C}_{\zeta},\\
\overline{\mat{\kappa}}_{n+1}^{\,-} &= \overline{\mat{\kappa}}_{n}^{\,+},\\
\mat{C}_{\kappa|n}^{\,-} &= \mat{C}_{\kappa|n}^{\,+},\\
\mat{C}_{n+1,\kappa|n}^{\, -} &= \overline{\mat{A}}_n^{\,+}\,\mat{C}_{n,\kappa|n}^{\, +} \;+\; \mat{J}_n^{\, +}\,\mat{C}_{\kappa|n}^{\,+}.
\end{align}
\end{subequations}
}
\Statex  \hspace{1.2em} \textbf{Update step:}
{
\setlength{\jot}{6pt}   
\begin{subequations}\label{eq:kalman-update}
\begin{align}
\overline{\ary{v}}_{n+1}^{\, +} 
&= \overline{\ary{v}}_{n+1}^{\, -} 
+ \mat{C}_{n+1}^{\, -}\,\mat{H}^\trans
\big(\mat{H}\,\mat{C}_{n+1}^{\, -}\,\mat{H}^\trans + \mat{C}_e\big)^{-1}
\big(\ary{y}_{n+1} - \mat{H}\,\overline{\ary{v}}_{n+1}^{\, -}\big),\\
\mat{C}_{n+1}^{\, +}             
&= \mat{C}_{n+1}^{\, -} 
- \mat{C}_{n+1}^{\, -}\,\mat{H}^\trans
\big(\mat{H}\,\mat{C}_{n+1}^{\, -}\,\mat{H}^\trans + \mat{C}_e\big)^{-1}
\mat{H}\,\mat{C}_{n+1}^{\, -},\\
\overline{\ary{\kappa}}_{n+1}^{\, +} 
&= \overline{\ary{\kappa}}_{n+1}^{\, -} 
+ \mat{C}_{\kappa,n+1|n}^{\, -}\,\mat{H}^\trans
\big(\mat{H}\,\mat{C}_{n+1}^{\, -}\,\mat{H}^\trans + \mat{C}_e\big)^{-1}
\big(\ary{y}_{n+1} - \mat{H}\,\overline{\ary{v}}_{n+1}^{\, -}\big),\\
\mat{C}_{\kappa|n+1}^{\, +}
&= \mat{C}_{\kappa|n+1}^{\, -}
- \mat{C}_{\kappa,n+1|n}^{\, -}\mat{H}^\trans
\big(\mat{H}\mat{C}_{n+1}^{\, -}\mat{H}^\trans + \mat{C}_e\big)^{-1}
\mat{H}\,\mat{C}_{n+1,\kappa|n}^{\, -},\\
\mat{C}_{n+1,\kappa|n}^{\, +}
&= \mat{C}_{n+1,\kappa|n}^{\, -}
- \mat{C}_{n+1}^{\, -}\mat{H}^\trans
\big(\mat{H}\mat{C}_{n+1}^{\, -}\mat{H}^\trans + \mat{C}_e\big)^{-1}
\mat{H}\,\mat{C}_{n+1,\kappa|n}^{\, -}.
\end{align}
\end{subequations}
}
\end{algorithmic}
\end{algorithm}

\subsection{Hyperparameter estimation \label{sec:bayes-params}}
%--------------------------------------------------------------------------------
The introduced statistical formulation depends on a set of hyperparameters denoted by $\boldsymbol{\theta}$. These include, for instance, parameters in the RF-SPDE formulation (Section~\ref{sec:forwardgenmat}) governing the mean and covariance of the forcing related.  The statistical model is calibrated by learning $\boldsymbol{\theta}$ from observation data. Depending on which intrinsic properties the hyperparameters represent for the statistical model, they can be assumed to be time-invariant $\boldsymbol{\theta}$ or time-variant $\boldsymbol{\theta}_n$ (but uncorrelated through time $p(\boldsymbol{\theta}_n|\boldsymbol{\theta}_{n-1})=p(\boldsymbol{\theta}_n)$).

When the hyperparameter dependence is made explicit, the posterior density~\eqref{eq:update-bayes} can be written as  
\begin{equation}\label{eq:update-bayes-theta}
p(\ary{v}_{n+1},\ary{\kappa}|\ary{y}_{1:n+1},\boldsymbol{\theta}_{{n}})
= \frac{p(\ary{y}_{n+1}|\ary{v}_{n+1},\ary{\kappa},\boldsymbol{\theta}_{{n}})\,p(\ary{v}_{n+1},\ary{\kappa}|\ary{y}_{1:n},\boldsymbol{\theta}_{{n}})}
       {p(\ary{y}_{n+1}|\ary{y}_{1:n},\boldsymbol{\theta}_{{n}})},
\end{equation}
representing the Level 1 state posterior of a hierarchical Bayes. For notational consistency with the graphical model (see Fig.~\ref{fig:bayesian_framework}), the hyperparameters estimated using data up to $t_{n+1}$ are denoted by $\boldsymbol{\theta}_{n}$.

In statFEM, we seek the posterior density by marginalising the factored joint between the hyperparameters and state
\begin{equation}
p(\ary{v}_{n+1},\ary{\kappa}|\ary{y}_{1:n+1})=  \int p(\ary{v}_{n+1},\ary{\kappa},\boldsymbol{\theta}_n|\ary{y}_{1:n+1}) \,\D \boldsymbol{\theta}_n
= \int p(\ary{v}_{n+1},\ary{\kappa}|\ary{y}_{1:n+1},\boldsymbol{\theta}_{{n}}) p(\boldsymbol{\theta}_{{n}}|\ary{y}_{1:n+1})\,\D \boldsymbol{\theta}_{{n}}.
\end{equation}
In the empirircal Bayes, or evidence approximation approach (see e.g.~\cite{mackay1999comparison}), this integral is approximated using a single evaluation point:
\begin{equation}
       p(\ary{v}_{n+1},\ary{\kappa}|\ary{y}_{1:n+1}) \approx p(\ary{v}_{n+1},\ary{\kappa}|\ary{y}_{1:n+1},\boldsymbol{\theta}^{*}_{n}),
\end{equation}
 A canonical choice for $\boldsymbol{\theta}^{*}_{n}$ {is the \textit{maximum a posteriori} (MAP) estimate}
\begin{equation}\label{eq:map}
\boldsymbol{\theta}^{*}_n = \arg\max_{\boldsymbol{\theta}_{{n}}}\, p(\boldsymbol{\theta}_{{n}}|\ary{y}_{1:n+1}), 
\end{equation}
%or the mean 
%\begin{equation}
%\boldsymbol{\theta}^{*}_{n} = \expect{[p(\theta|\ary{y}_{1:n+1})]}. 
%\end{equation}
where the Level 2 hyperparameter posterior $p(\boldsymbol{\theta}_n|\ary{y}_{1:n+1})$ is maximized up to a proportionality constant for the MAP estimate:
\begin{equation}
       p(\boldsymbol{\theta}_{{n}}|\ary{y}_{1:n+1}) =  
       \frac{p(\ary{y}_{n+1}|\ary{y}_{1:n},\boldsymbol{\theta}_{n})p(\boldsymbol{\theta}_{n})}
       {\int p(\ary{y}_{n+1}|\ary{y}_{1:n},\boldsymbol{\theta}_{n}) p(\boldsymbol{\theta}_{{n}}|\ary{y}_{1:n})\, \D \boldsymbol{\theta}_{n}} \propto p(\ary{y}_{n+1}|\ary{y}_{1:n},\boldsymbol{\theta}_{{n}}) p(\boldsymbol{\theta}_{{n}}).
\end{equation}
The first term on the right-hand side is the marginal likelihood of the state, i.e. the expression in the denominator of the Bayes formula~\eqref{eq:update-bayes-theta}. The second term is the prior on the hyperparameters (or hyperprior). The marginal likelihood is also a Gaussian density for the considered linear observation model~\eqref{eq:linear-sys}, and can be obtained by marginalising the auguented state (cf. also~\eqref{eq:marg-like-bayes}):
\begin{equation}\label{eq:pred-like-new}
	\begin{aligned}
p\!\left(\ary{y}_{n+1}| \ary{y}_{1:n},\boldsymbol{\theta}_{{n}}\right)
&= \int p(\ary{y}_{n+1}|\ary{v}_{n+1},\ary{\kappa},\ary{\theta}_{n})\,p(\ary{v}_{n+1},\ary{\kappa}|\ary{y}_{1:n},\ary{\theta}_{n})\,\D\ary{v}_{n+1}\,\D\ary{\kappa}\\
&=~
\mathcal{N}\,\Big(
\mat{H}\,\overline{\ary{v}}_{n+1}^{\, -}(\boldsymbol{\theta}_{{n}}),
~\mat{H}\mat{C}_{n}^{\, -}(\boldsymbol{\theta}_{{n}})\mat{H}^\trans
   + \mat{C}_e(\boldsymbol{\theta}_{{n}})
\Big).
	\end{aligned}
\end{equation}

In our numerical examples, we optimise for $\boldsymbol{\theta}^{*}_{n}$  equivalently by minimising the negative logarithm of~\eqref{eq:map}:
\begin{equation}
\phi_{n+1}(\boldsymbol{\theta}_{{n}})
~=~
-\log p\!\left(\ary{y}_{n+1}| \ary{y}_{1:n},\boldsymbol{\theta}_{{n}}\right)
-\log p(\boldsymbol{\theta}_{{n}}).
\end{equation}
Taking the logarithm of~\eqref{eq:pred-like-new}, the negative marginal log-likelihood is
\begin{align}\label{eq:neg-log-like}
-\log p\!\left(\ary{y}_{n+1}| \ary{y}_{1:n},\boldsymbol{\theta}_{{n}}\right)
&=~\tfrac{1}{2}
\Big[
\big(\ary{y}_{n+1}-\mat{H}\,\overline{\ary{v}}_{n+1}^{\, -}(\boldsymbol{\theta}_{{n}})\big)^\trans
\!\big(
\mat{H}\mat{C}_{n+1}^{\, -}(\boldsymbol{\theta}_{{n}})\mat{H}^\trans
+\mat{C}_e(\boldsymbol{\theta}_{{n}})
\big)^{-1}\!
\big(\ary{y}_{n+1}-\mat{H}\,\overline{\ary{v}}_{n+1}^{\, -}(\boldsymbol{\theta}_{{n}})\big) \notag\\
&\quad
+ \log\!\left|\big(
\mat{H}\mat{C}_{n+1}^{\, -}(\boldsymbol{\theta}_{{n}})\mat{H}^\trans
+\mat{C}_e(\boldsymbol{\theta}_{{n}})
\big)\right|
+ n_y\log(2\pi)
\Big],
\end{align}
where $n_y$ is the dimension of $\ary{y}_{n+1}$.

If $\boldsymbol{\theta}$ is assumed constant in time, the marginal likelihood over all observations is given by
\begin{equation}\label{eq:static-post}
p(\boldsymbol{\theta}| \ary{y}_{1:N})
~\propto~
p(\boldsymbol{\theta})
\prod_{n=1}^{N}
p\!\left(\ary{y}_n| \ary{y}_{1:n-1},\boldsymbol{\theta}\right),
\end{equation}
which follows directly from the probability chain rule. Taking the negative logarithm yields
\begin{equation}\label{eq:phi-batch}
\phi_{1:N}(\boldsymbol{\theta})
~=~
-\!\sum_{n=1}^{N}
\log p\!\left(\ary{y}_n| \ary{y}_{1:n-1},\boldsymbol{\theta}\right)
-\log p(\boldsymbol{\theta}).
\end{equation}
When estimating $\boldsymbol{\theta}$, the marginal likelihood is based on an augmented state; thus, $\ary{\kappa}$ is jointly updated with $\boldsymbol{\theta}$ from observations as well.

It is noted that the filtering procedure with an augmented material state provides a joint parameter estimation framework for the material field parameters $ \overline{\ary{\kappa}} $ and $ \mat{C}_\kappa $, which are updated recursively as new observations are assimilated. Alternatively, these parameters (i.e. hyperparameters in this case) can be obtained via MAP estimation by setting $\{\boldsymbol{\theta},\ary{\kappa}\}\to \boldsymbol{\theta}$. The state augmentation approach is adopted here as it yields an explicit joint posterior $ p(\ary{v}_n, \ary{\kappa} | \vec{y}_{1:n}) $, enabling coherent Bayesian information flow between the state and material fields while avoiding computationally demanding non-convex optimisation. We note that, due to zero-process-noise dynamics on $\ary{\kappa}$, the posterior variance may degenerate ($ \mat{C}_{\kappa|n}^{+} \to \mat{0} $) \cite{sarkka2023bayesian}. In the idealised limit of informative observations, the posterior mean converges to the true value, $ \overline{\ary{\kappa}}_n^{+} \to \ary{\kappa}_{\text{true}} $.

\subsection{Parameter estimation in the SDOF system}
\begin{figure}[!t]
    \centering
    \subfloat[$\sigma_f$ estimation \label{fig:likelihood}]
        {\includegraphics[width=0.48\textwidth]{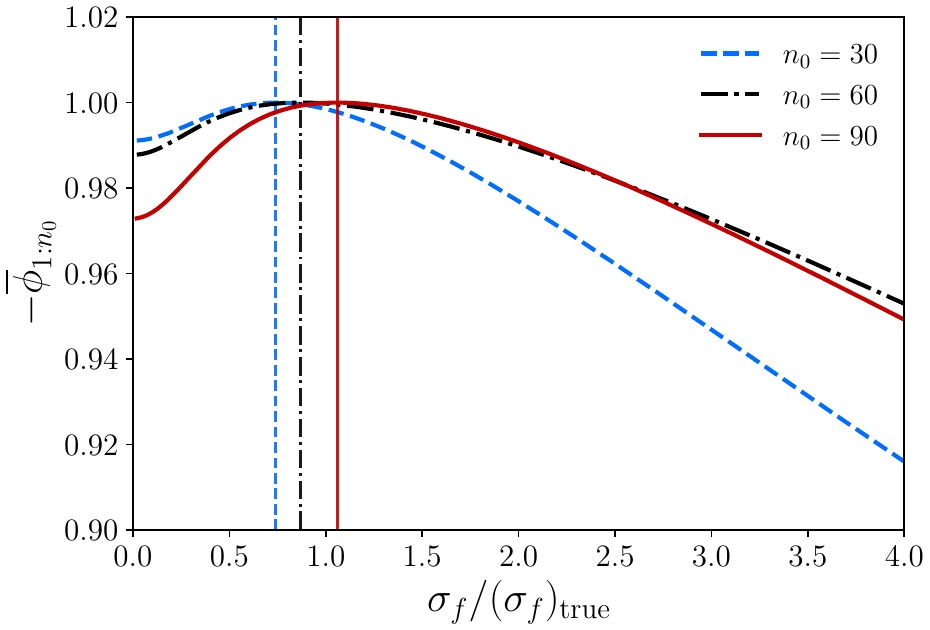}\vspace{-1mm}}
    \hfill
    \subfloat[Stiffness update ($n_0=90$) \label{fig:stiff_up}]
        {\includegraphics[width=0.48\textwidth]{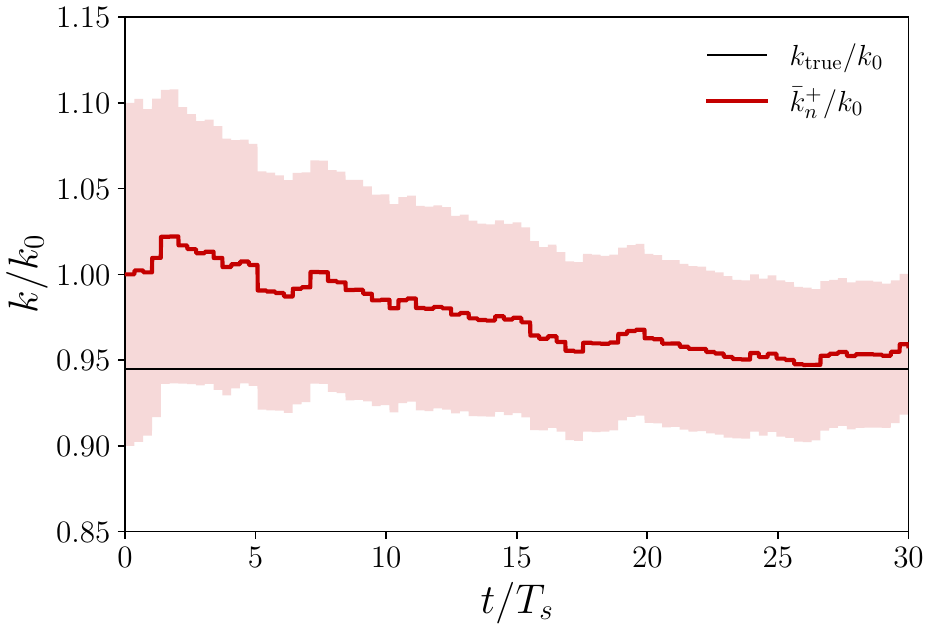}\vspace{-1mm}}
    \caption{(a) Normalized marginal negative log likelihood $-\phi_{1:n_0}(\sigma_f)$ for different numbers of observations $n_0$, shown as a function of the forcing variance $\sigma_f$. The vertical lines correspond to the estimated $\hat{\sigma}_f$ where the likelihood is maximal. (b) The updated stiffness value $\overline{k}_n^{+}$ as the observation data is read along with its updated confidence interval. }
    \label{fig:stiff_likelihood}
\end{figure}

We use the SDOF example from Section~\ref{sec:forward} to demonstrate the estimation of the misspecification standard deviation $\sigma_f$, treated here as a time-invariant hyperparameter endowed with a uniform prior. The observation data are generated using a true stiffness value of $k_{\mathrm{true}}=94.48$ and true forcing variance value of $(\sigma_f)_{\mathrm{true}} = 0.05$ with $\sigma_e=0.005$. Figure~\ref{fig:likelihood} shows the resulting marginal likelihoods for two measurement-noise levels $\sigma_e$ and three observation sampling factors $n_o = \Delta t / \Delta t_o$, where $\Delta t$ and $\Delta t_o$ denote the simulation and observation time steps. The shape of the likelihood and the location of its maximizer (the MAP estimate $\hat{\sigma}_f$) depend on the sampling factor $n_o$ and, as expected, increasing the number of observations improves the accuracy of $\hat{\sigma}_f$. Figure~\ref{fig:stiff_up} shows the updated $\overline{k}_n^{+}$ obtained using the MAP estimation of $\sigma_f$, starting with the prior mean value $k_0=100$ at $t=0$, and converging towards the $k_{\mathrm{true}}$ as more observation data is read over time. Additionally, we also notice that the updated variance, $(\sigma_{k})_n^+$, is reducing as is evident from the confidence intervals. 

%% file: examples.tex
%
%--------------------------------------------------------------------------------          
\section{Examples \label{sec:examples}}
%--------------------------------------------------------------------------------
%
In this section, we demonstrate the proposed Bayesian filtering framework through two examples. In Section~\ref{sec:examplesbar}, we consider a one-dimensional elastodynamic problem with uncertain material properties and stochastic forcing, and examine how the material correlation length influences state estimation and parameter updating, including a comparison with a filtering approach in which the material field is not updated. In Section~\ref{sec:exampleshole}, we study a two-dimensional anti-plane elasticity problem with similar sources of uncertainty, and assess the ability of the framework to recover displacements, including at unobserved locations, along with the associated material field.
%
%--------------------------------------------------------------------------------          
\subsection{Bar with axial harmonic loading \label{sec:examplesbar}}
%--------------------------------------------------------------------------------
%
\subsubsection{Problem setup}
\label{sec:ps1}

\begin{figure}[htbp]
	\centering
	\includegraphics[width=0.65\textwidth]{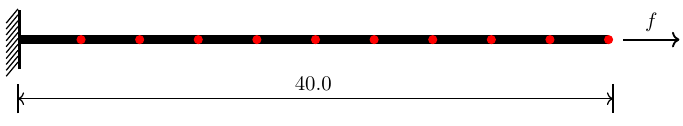}
	\caption{Schematic of a bar with axial harmonic loading. The red dots indicate the locations of sensor measurements.}
	\label{fig:Ex1_Schematic_1DBar}
\end{figure}

We consider a fixed-end bar of density, $\rho=1200$, and of length $L=40.0$ discretized into $n_{\mathrm{el}}=80$ elements with uniform element size $h=L/n_{\mathrm{el}}=0.5$ (see Fig.~\ref{fig:Ex1_Schematic_1DBar}). The bar is clamped at the left end, and the axial load is applied at the tip (right end). The mean forcing is sinusoidal,
\begin{equation}
\overline{f}(t) = f_0\sin(2\pi f_b t),
\qquad
f_0=2\times 10^3,
\qquad
f_b=0.25,    
\end{equation}
with $f_b$ lying between the first two axial-mode frequencies, $f_1=0.12$ and $f_2=0.38$. In addition, the tip load is perturbed by an additive white-noise component with standard deviation $\sigma_f$. We treat $\sigma_f$ as a hyperparameter, with the true value taken as $\sigma_f = 0.05f_0$. Damping is modeled using Rayleigh damping calibrated from the two target modal frequencies $f_1$ and $f_2$ with a damping ratio of $0.5\%$. 

Spatial variability in stiffness is introduced through a zero-mean Gaussian random field $\kappa(x)$. The Young’s modulus is then defined through the transformation
\begin{equation}
E = \widetilde{E}\,\exp(\kappa),
\qquad
\kappa \sim \mathcal{N}(0,\,c_\kappa),    
\end{equation}
which ensures that it remains positive, where $c_\kappa$ is constructed from an SPDE-based Mat\'ern covariance with smoothness parameter $\nu=1.5$, marginal standard deviation $\sigma_\kappa=0.1$, and correlation length $l_\kappa$ (see~\ref{sec:matern}). To ensure that the prior mean of $E$ equals a prescribed value ${E}_0$, we set
\begin{equation}
\widetilde{E} = \frac{{E}_0}{\exp\!\left(\tfrac{1}{2}\sigma_\kappa^2\right)},
\qquad
E_0=5\times 10^5,    
\end{equation}
so that $\mathbb{E}[E]=E_0$. Let the the longitudinal wave speed in the nominal material be denoted as $v_l = \sqrt{{E_0}/{\rho}}$, nominal time period as $T_s=L/v_l$ and the nominal displacement as $u_s = f_0 L/E_0$. 

The time integration uses a constant time step $\Delta t/T_s = 2.5\times10^{-3}$ over a total simulation time $T=10T_s$ with $u=0$ and $\dot{u}=0$ throughout the domain as initial conditions. Observations are sampled every $n_o=1000$ time steps, corresponding to an observation interval $\Delta t_o = n_o \Delta t = 0.5$. We use $n_s=10$ displacement sensors placed approximately uniformly along the bar (see Fig.~\ref{fig:Ex1_Schematic_1DBar}), and observations are assumed to be affected by independent Gaussian noise with standard deviation $\sigma_e = 10^{-2}$. We investigate the effect of correlation length by considering several values of $l_\kappa$. Larger $l_\kappa$ induces smoother stiffness variations, which improves identifiability of the material property field from sparse observations. The reported results compare: (i) the posterior predictive displacement time histories at the tip, for different choices of the correlation length $l_{\kappa}$ and (ii) the posterior estimate of the Young’s modulus field $E$. To avoid strong conditioning at early times, observations are introduced only after an initial transient period during which the system evolves without measurement updates. Additionally, measurements are not taken after a certain time to study how the state evolves after certain number of observations. 

\subsubsection{Parameter estimation and filtering results}
\label{sec:spe1}
\begin{figure}[!t]
    \centering
	\subfloat[$l_\kappa=2.5$\, \label{fig:forward_1d_l2.5}]
        {\includegraphics[width=0.48\textwidth]{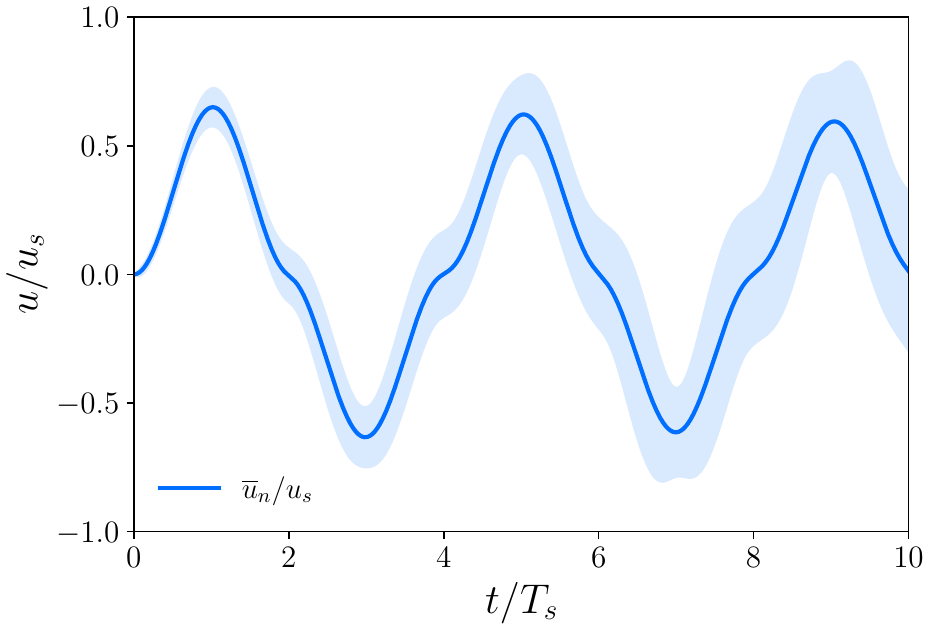}\vspace{-1mm}}\hfill
	\subfloat[$l_\kappa=10.0$\, \label{fig:forward_1d_l10}]
        {\includegraphics[width=0.48\textwidth]{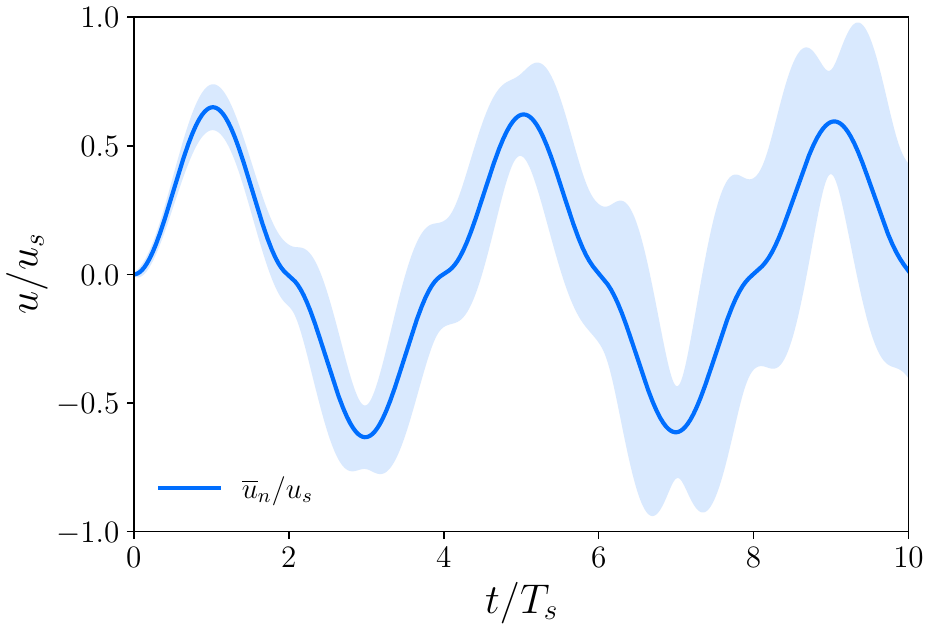}\vspace{-1mm}}
	\caption{Forward predictive displacement at the tip of the bar for two correlation lengths of the material field.}
	\label{fig:forward_1d}
\end{figure}

\begin{figure}[!t]
    \centering
	\subfloat[$l_\kappa=2.5$\, \label{fig:likelihood_l2.5}]
        {\includegraphics[width=0.48\textwidth]{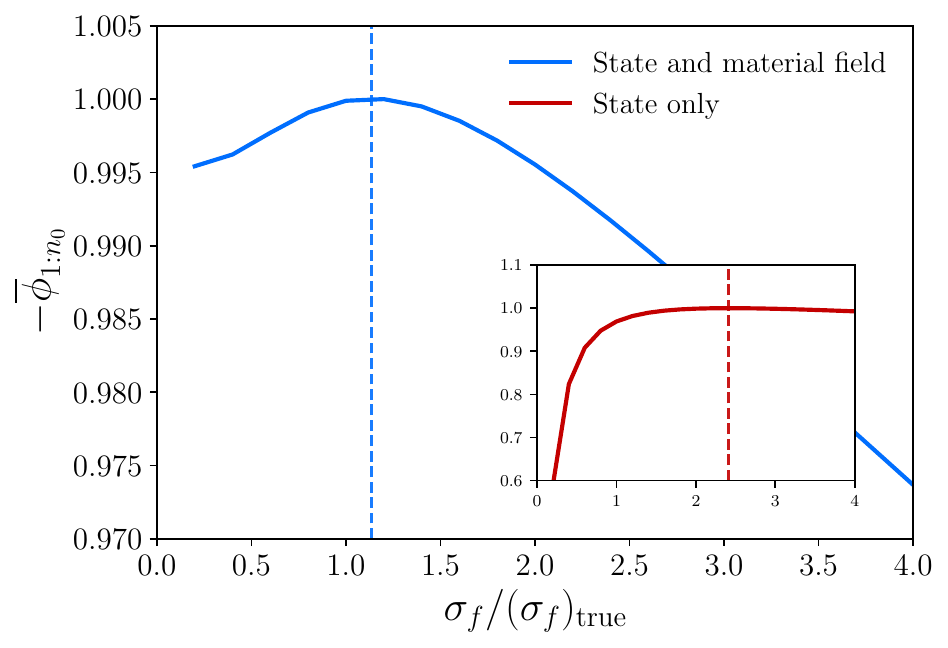}\vspace{-1mm}}\hfill
	\subfloat[$l_\kappa=10.0$\, \label{fig:likelihood_l}]
        {\includegraphics[width=0.48\textwidth]{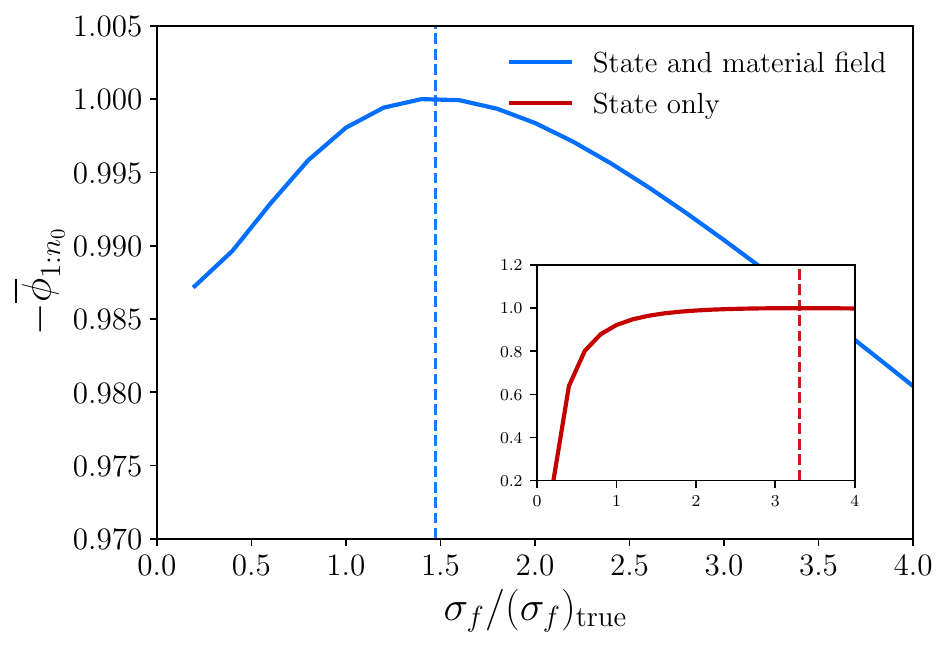}\vspace{-1mm}}
	\caption{Marginal likelihood as a function of the misspecification standard deviation $\sigma_f$ for two correlation lengths of the material field for the augmented filtering formulation. The insets show the likelihood when the material field is not updated.}
	\label{fig:likelihood_1d}
\end{figure}

\begin{figure}[!t]
    \centering
	\subfloat[Updated displacement ($l_\kappa=2.5$) \label{fig:update_disp_1d_l2.5}]
        {\includegraphics[width=0.48\textwidth]{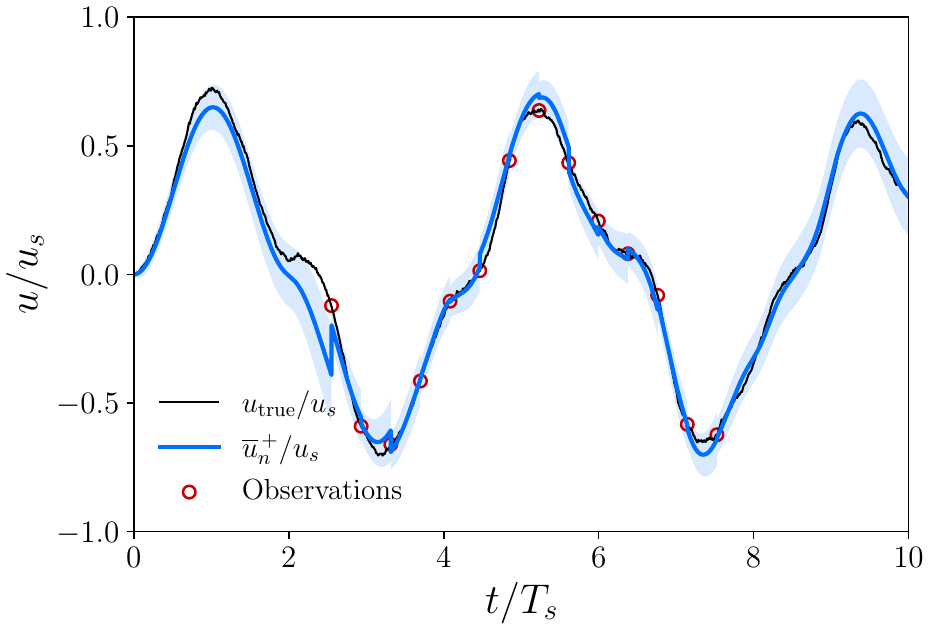}\vspace{-1mm}}\vfill
    \subfloat[Updated Young’s modulus field at $t=0$ and $t=5T_s$ ($l_\kappa=2.5$) \label{fig:update_E_1d_his_l2.5}]
        {\includegraphics[width=0.48\textwidth]{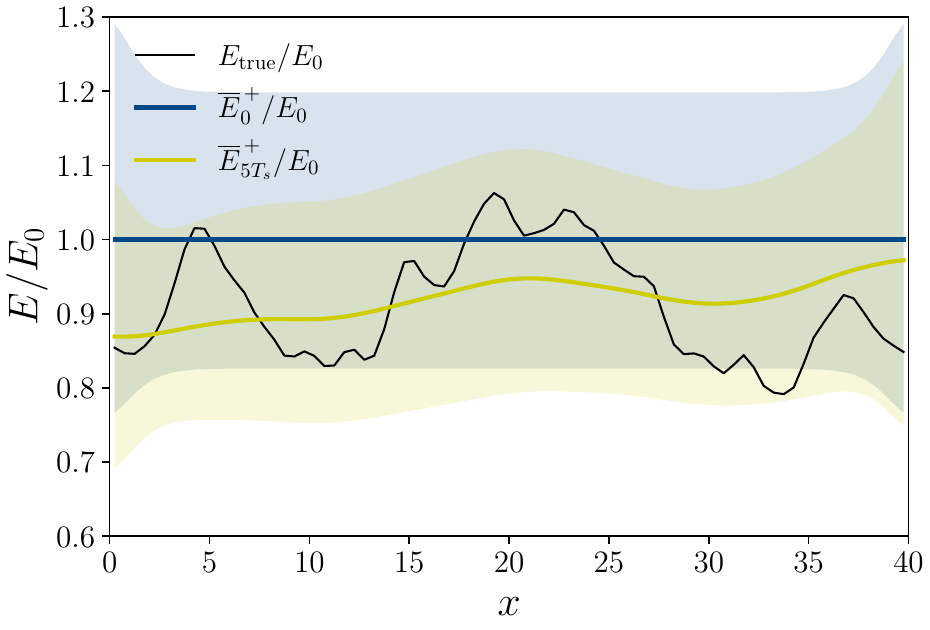}\vspace{-1mm}}
    \hfill        
    \subfloat[Updated Young’s modulus field at $t=10T_s$ ($l_\kappa=2.5$) \label{fig:update_E_1d_l2.5}]
        {\includegraphics[width=0.48\textwidth]{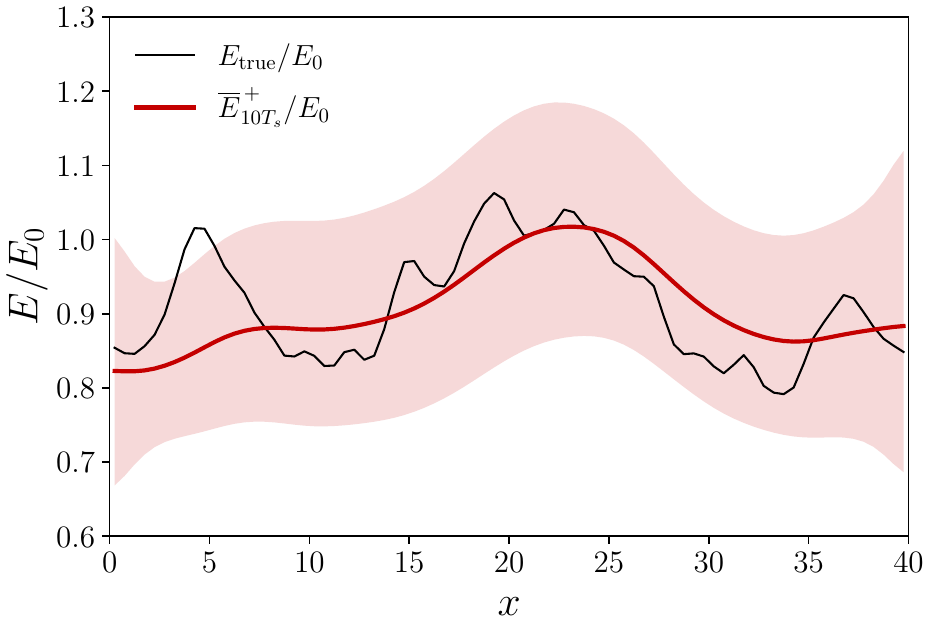}\vspace{-1mm}}
        \caption{Posterior predictive displacement and updated Young’s modulus field obtained from the augmented state and parameter estimation framework ($l_\kappa=2.5$).}
	\label{fig:update_1d_l2.5}
\end{figure}

\begin{figure}[!t]
    \centering
	\subfloat[Updated displacement ($l_\kappa=10.0$) \label{fig:update_disp_1d_l10}]
        {\includegraphics[width=0.48\textwidth]{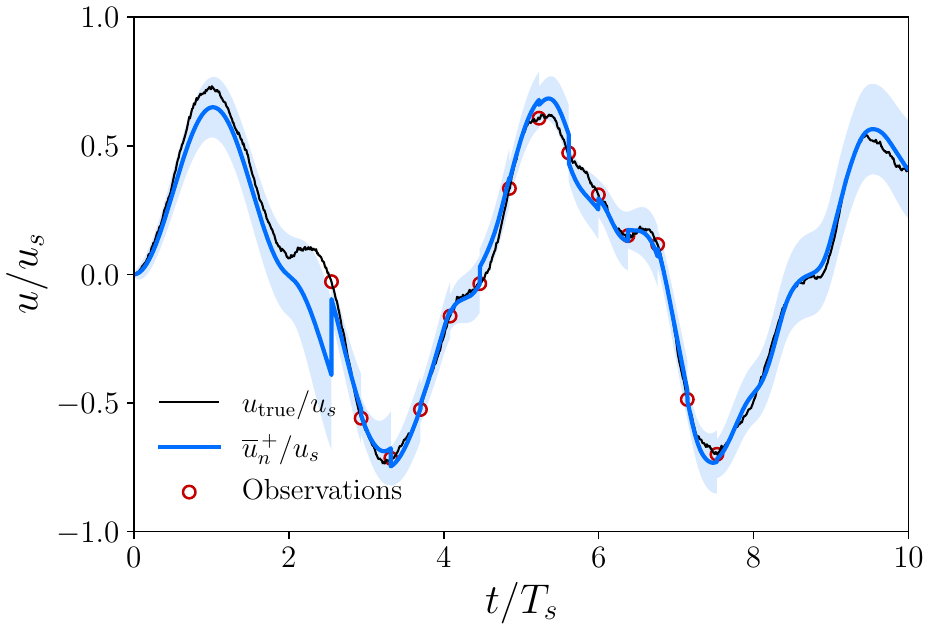}\vspace{-1mm}}
	\vfill	
	\subfloat[Updated Young’s modulus field at $t=0$ and $t=5T_s$ ($l_\kappa=10.0$) \label{fig:update_E_1d_his_l10}]
        {\includegraphics[width=0.48\textwidth]{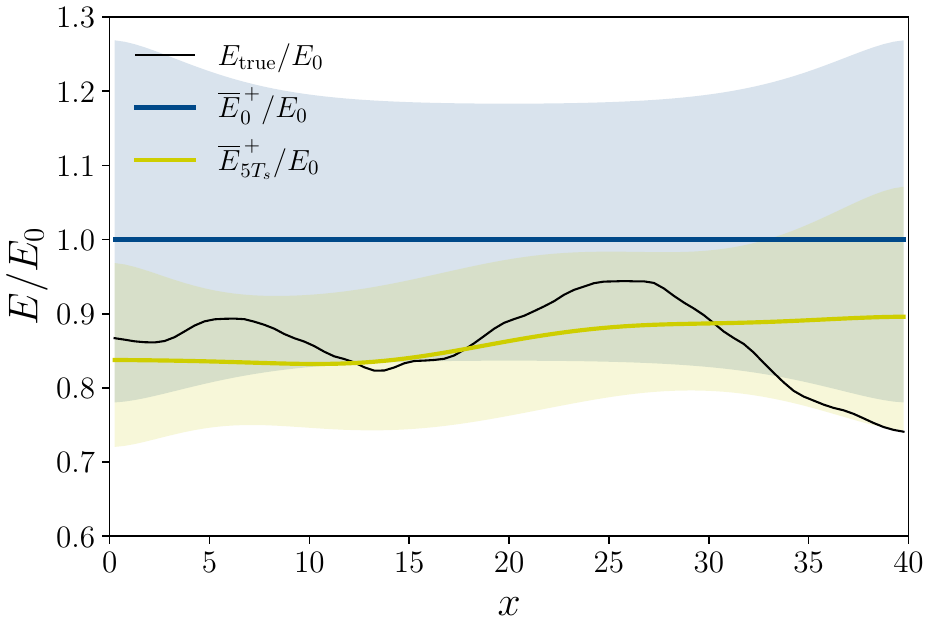}\vspace{-1mm}}\hfill
    \subfloat[Updated Young’s modulus field at $t=10T_s$ ($l_\kappa=10.0$) \label{fig:update_E_1d_l10}]
        {\includegraphics[width=0.48\textwidth]{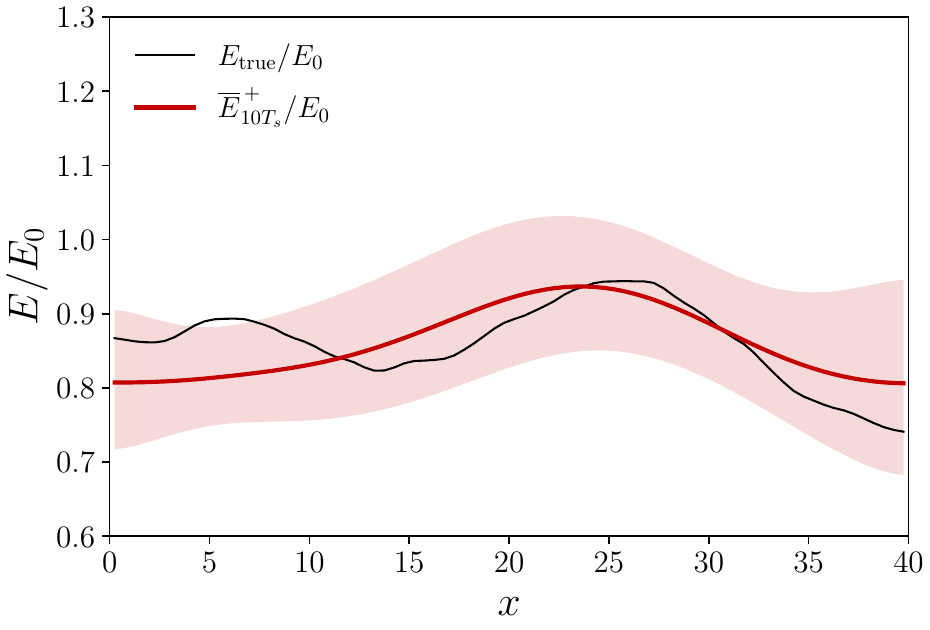}\vspace{-1mm}}	
	\caption{Posterior predictive displacement and updated Young’s modulus field obtained from the augmented state and parameter estimation framework ($l_\kappa=10.0$).}
	\label{fig:update_1d_l10}
\end{figure}

The parameter estimation and filtering results for the two choices of the correlation length $l_\kappa$ of the material field $\kappa$ are summarized in Figures~\ref{fig:forward_1d}--\ref{fig:update_1d_l10}. Figure~\ref{fig:forward_1d} shows the forward predictive displacement at the tip of the bar obtained from the stochastic model~\eqref{eq:mean-update} and~\eqref{eq:cov-compact}, incorporating both material and forcing uncertainty. A clear dependence on the correlation length is observed. For the shorter correlation length $l_\kappa=2.5$, the predictive variance of the tip displacement is smaller than for the longer correlation length $l_\kappa=10.0$. This behavior can be attributed to the nature of the underlying random field realizations. Shorter correlation lengths for $\kappa$ lead to rapidly fluctuating stiffness fields whose effects tend to average out over the domain. In contrast, longer correlation lengths induce coherent spatial trends in $\kappa$, which act as effective stiffness variations over large portions of the domain and allow the system response to deviate more strongly from its nominal mean. 

Figure~\ref{fig:likelihood_1d} presents the estimation of the misspecification standard deviation $\sigma_f$ for the two cases of different correlation lengths considered using the marginal likelihood approach described in Section~\ref{sec:bayes-params}. In both cases, the likelihood exhibits a well-defined maximum close to the true value of $\sigma_f$, indicating that the misspecification standard deviation is accurately identified from the displacement observations. Additionally, the inset plots show the estimation of $\sigma_f$ for the case when the material property field is not updated along with the displacement and velocity fields. In contrast to the augmented case, the marginal likelihood is no longer sharply peaked and its maximum is shifted away from the true value of $\sigma_f$.

Figures~\ref{fig:update_1d_l2.5} and~\ref{fig:update_1d_l10} illustrate the updated displacement and material property field results obtained using the estimated misspecification variance values within the augmented filtering framework~\eqref{eq:kalman-pred} and~\eqref{eq:kalman-update}. The posterior predictive displacement at the tip is shown together with the updated Young’s modulus field at different time instants for both correlation lengths. In both cases, the posterior mean response deviates from the forward-model prediction initially and closely matches the true displacement once observations are assimilated. Concurrently, the posterior uncertainty in both the displacement and the material property fields decreases as additional measurements are incorporated. Once the observations are not read anymore, the posterior uncertainty in displacement increases. A notable distinction between the two cases arises in the inferred Young’s modulus field. For the shorter correlation length $l_\kappa=2.5$, the posterior variance of the modulus remains comparatively larger, reflecting the difficulty of inferring rapidly fluctuating material fields from sparse observations. Conversely, for $l_\kappa=10.0$, the stronger spatial correlation leads to smoother posterior fields with reduced uncertainty, indicating improved identifiability. 

\begin{figure}[!t]
    \centering
	\subfloat[Updated displacement ($l_\kappa=2.5$) \label{fig:update_disp_1d_so_l2.5}]
        {\includegraphics[width=0.48\textwidth]{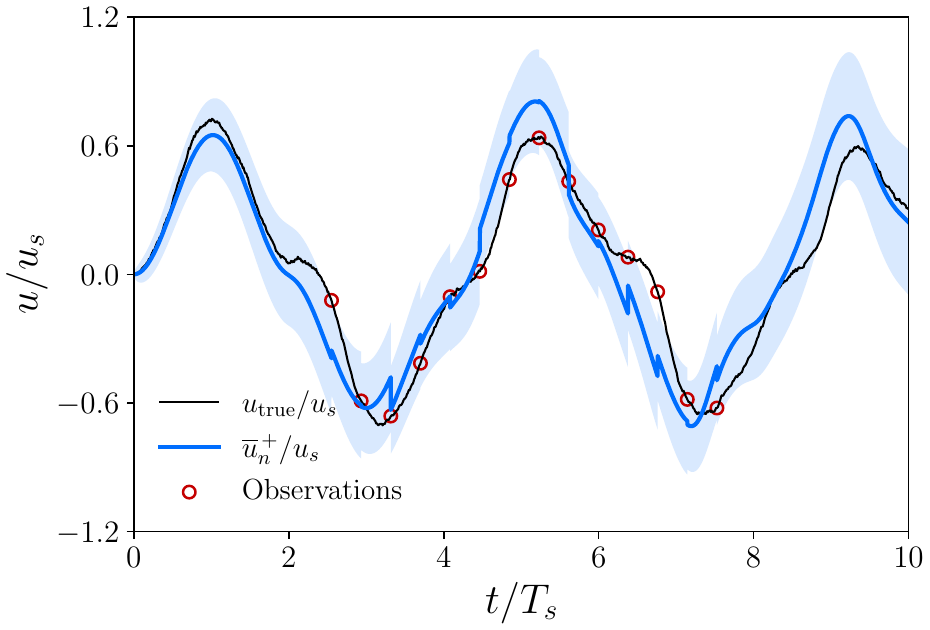}\vspace{-1mm}}
	\hfill	
	\subfloat[Updated displacement ($l_\kappa=10.0$) \label{fig:update_disp_1d_so_l10}]
        {\includegraphics[width=0.48\textwidth]{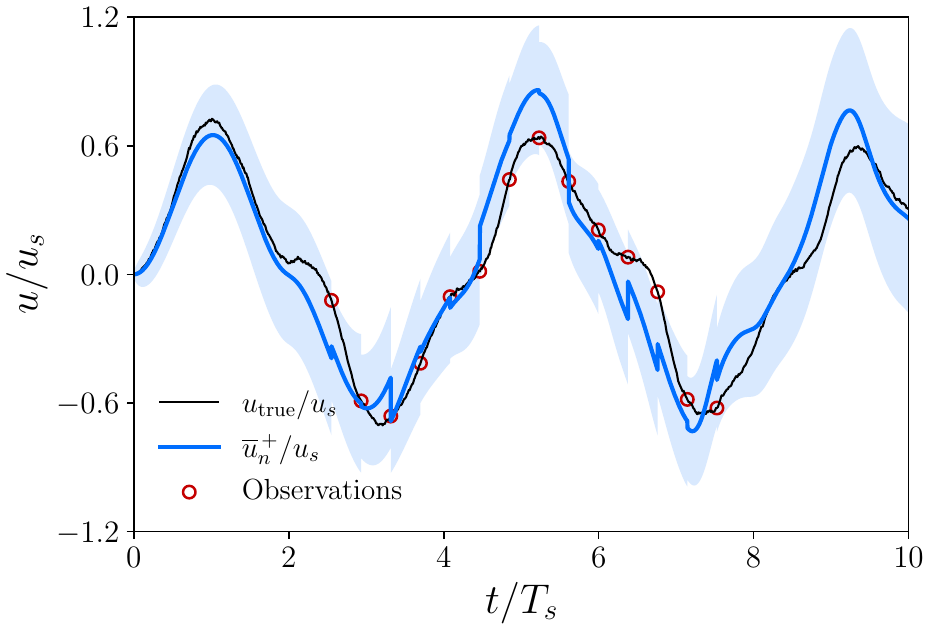}\vspace{-1mm}}	
	\caption{Posterior predictive displacement without updating the Young's modulus field.}
	\label{fig:update_1d_so}
\end{figure}

Figure~\ref{fig:update_1d_so} examines the effect of assimilating displacement and velocity observations while the material field is not updated; the misspecification standard deviation $\hat{\sigma}_f$ is taken as the value estimated from the non-augmented filtering formulation. In contrast to Figs.~\ref{fig:update_1d_l2.5} and~\ref{fig:update_1d_l10}, where the filter jointly updates both the displacement and the Young’s modulus, the posterior displacement obtained without updating $\kappa$ shows reduced agreement with the true response, with discrepancies becoming more pronounced for larger correlation lengths. This reflects the fact that, for longer correlation lengths, the system response is strongly influenced by spatially coherent variations in material properties. When $\kappa$ is not updated, these effects cannot be represented through the material field. As a result, the model--data discrepancy is absorbed into the state estimate rather than attributed to spatial variations in material properties, leading to a less accurate reconstruction of the response. This comparison highlights the importance of joint state and material property updating, particularly when material heterogeneity exhibits strong spatial correlation.
%
%--------------------------------------------------------------------------------          
\subsection{Anti-plane elasticity of a plate with hole \label{sec:exampleshole}}
%--------------------------------------------------------------------------------
\subsubsection{Problem setup}
\label{sec:ps2}
\begin{figure}[!t]
    \centering
    \subfloat[Geometry and boundary conditions\label{fig:dom_anti}]
    {\includegraphics[width=0.45\textwidth]{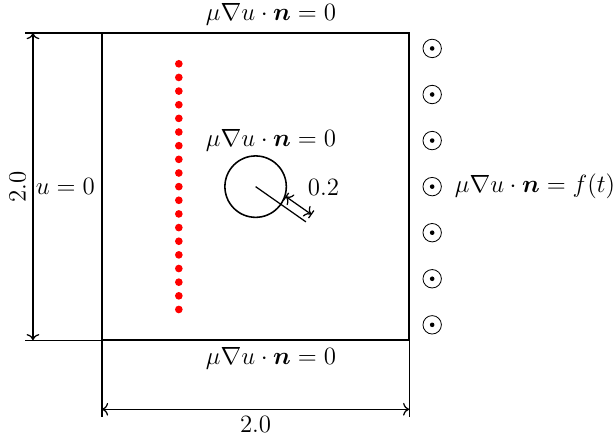}\vspace{-1mm}}
    \hfill
    \subfloat[Triangular pulse $\overline{f}(t)$ used for the deterministic traction \label{fig:pulse_anti}]
    {\includegraphics[width=0.55\textwidth]{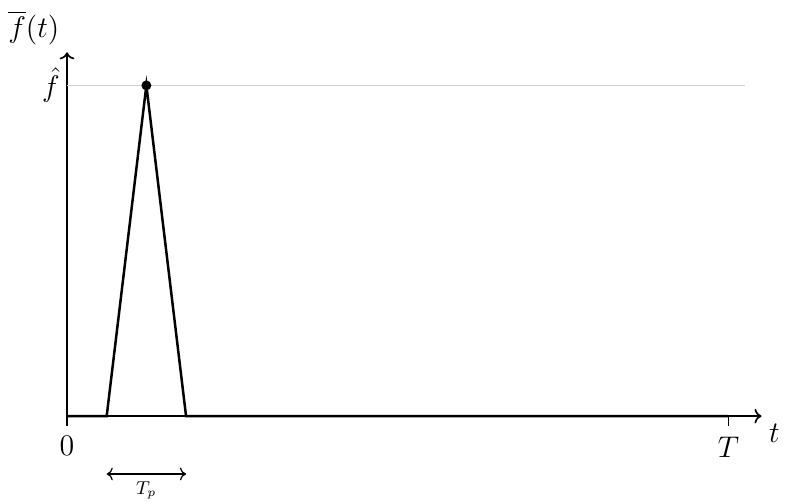}\vspace{-1mm}}
    \vfill
    \subfloat[True shear modulus field $\mu(\vec{x})$ \label{fig:shear_field}]
    {\includegraphics[width=0.5\textwidth]{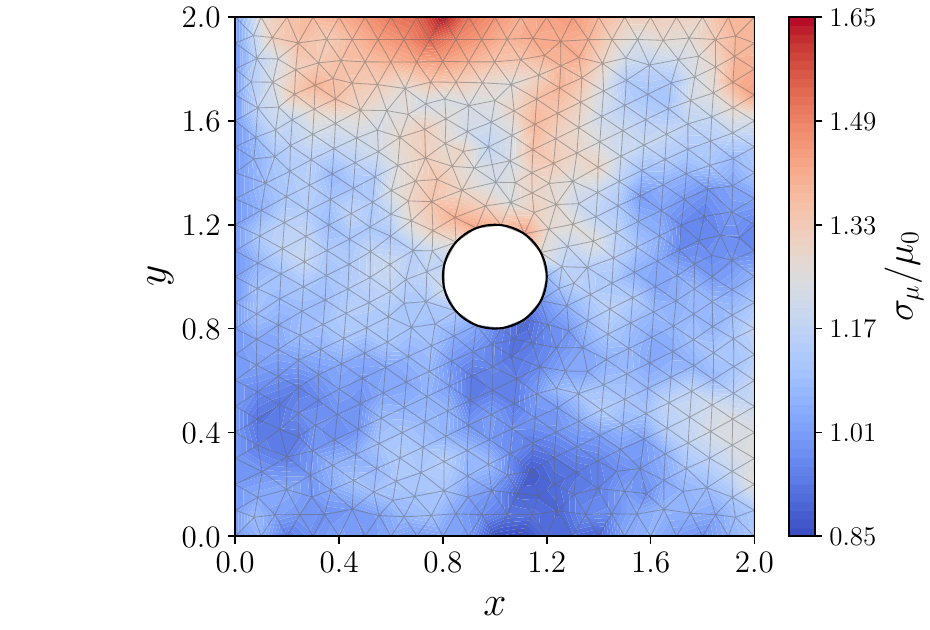}\vspace{-1mm}}
    \caption{Schematic of the anti-plane shear problem: (a) geometry and boundary conditions; the red dots indicate the locations of sensor measurements, (b) deterministic traction pulse and (c) true shear modulus field.}
    \label{fig:Schematic_anti}
\end{figure}

We consider a two-dimensional anti-plane shear problem defined on a square domain $\Omega = [0,L_x]\times[0,L_y]\setminus B_R,$ where $L_x = L_y = 2.0$ and $B_R$ denotes a circular hole of radius $R = 0.20$ centered at $(c_x,c_y) = (1.0,1.0)$. The anti-plane (out-of-plane) displacement field is denoted by $u(\vec{x},t)$, where $\vec{x}\in\Omega \subset \mathbb{R}^2$ and $t\in[0,T]$. The anti-plane shear modulus is denoted by $\mu(\vec{x})$ and we consider the density to be $\rho=8000$. The boundary on which Dirichlet conditions are applied corresponds to the left edge of the domain, while Neumann conditions are applied on the remaining boundaries (see Fig.~\ref{fig:dom_anti}).

\begin{comment}
The governing equations of elastodynamics~\eqref{eq:goveq} for anti-plane elasticity reduce to 
\begin{subequations}
	\begin{align}
    \nabla\cdot\!\big(\mu(\vec{x})\,\nabla u(\vec{x},t)\big)
	&= \rho\,\ddot{u}(\vec{x},t)
	\quad && \text{in } \Omega\times(0,T],	\\	
    u(\vec{x},t)&=0\quad &&\text{on }\{ (0,x_2)\} \times [0,T],\\
	\mu(\vec{x})\,\nabla u(\vec{x},t)\cdot\vec{n} &= f(t) \quad && \text{on }\{ (L_x,x_2)\} \times [0,T],\\
	\mu(\vec{x})\,\nabla u(\vec{x},t)\cdot\vec{n} &= 0 \quad && \text{on }\{ (x_1,0), (x_1,L_y)\}\cup B_R \times [0,T],\\
    u(\vec{x},0)=0\text{, } \dot{u}(\vec{x},0)&=0 \quad &&\text{for } \vec{x}\in\Omega,
\end{align}
\end{subequations}
where $\mu(\mathbf{x})$ is the shear modulus and we take $\rho = 8000$. 
\end{comment}
\newpage
The deterministic component of the applied traction, $\overline{f}(t)$, is taken as a triangular pulse (see Fig.~\ref{fig:pulse_anti}) with peak amplitude $\hat{f} = 5\times10^{5}$, starting at $t_0 = 0$,
\begin{subequations}
\begin{align}
\overline{f}(t) & =
\begin{cases}
\dfrac{\hat{f}}{T_r}(t - t_0),
& t_0 \le t < t_0 + T_r, \\[6pt]
\hat{f}\,\left(1 - \dfrac{t - (t_0 + T_r)}{T_f}\right),
& t_0 + T_r \le t \le t_0 + T_r + T_f, \\[6pt]
0,
& \text{otherwise}.
\end{cases}	\\
{f}_0 &= \frac{1}{T}\int_0^T \overline{f}(t)\,\mathrm{d}t
= \frac{(T_r+T_f)\,\hat{f}}{2T}.	
\end{align}
\end{subequations}

Spatial variability in the shear modulus is modeled using a lognormal random field,
\begin{equation}
\mu(\mathbf{x}) = \widetilde{\mu}\exp\!\big(\kappa(\vec{x})\big),	\qquad
\kappa \sim \mathcal{N}(0,\,c_\kappa), 
\end{equation}
where $c_\kappa$ is constructed from an SPDE-based Mat\'ern covariance with smoothness parameter $\nu=1.0$, standard deviation $\sigma_\kappa=0.1$, and correlation length $l_\kappa=1.0$ (see~\ref{sec:matern}). To ensure that the prior mean of $\mu(\vec{x})$ equals a prescribed value ${\mu}_0$, we set
\begin{equation}
\widetilde{\mu} = \frac{\mu_0}{\exp\!\left(\tfrac{1}{2}\sigma_\kappa^2\right)},
\qquad
\mu_0=2\times10^9,    
\end{equation}
so that $\mathbb{E}[\mu(\vec{x})]=\mu_0$. 

Let the the shear wave speed in the nominal material be denoted as $v_s = \sqrt{{\mu_0}/{\rho}}$, nominal time period as $T_s=L_y/v_s$, the nominal displacement as $u_s = f_0 L_y/\mu_0$. We set the rise and fall times of the pulse as $T_r = T_f = 0.5T_s$. Damping is modeled using Rayleigh damping, calibrated to achieve a damping ratio of $0.5\%$ at two target modal frequencies, $\omega_1 = v_s / L_x$ and $\omega_2 = v_s / (2h)$, where $h$ denotes the smallest mesh size. The stochastic component of the applied traction is a Gaussian white noise field, which is spatially correlated and temporally uncorrelated as described in sections~\ref{sec:forwardgoveq} and \ref{sec:forwardtime}. We chose the spatial correlation to be of Matérn type i.e. $c_f(\vec{x},\vec{x}^\prime) = \sigma_f^2 c_{\mathcal{M}}(\vec{x},\vec{x}^\prime)$ (see~\eqref{eq:Matern}). The correlation length is set as $l_f = 0.25L_y = 0.5$, the smoothness parameter as $\nu_f = 1.5$ and the marginal standard deviation $\sigma_f=0.05f_0$. 

The time integration uses a constant time step $\Delta t/T_s = 2\times10^{-2}$ over a total simulation time $T = 10T_s$, with $u = 0$ and $\dot{u} = 0$ throughout the domain as initial conditions. Displacement observations are collected at fixed interior sensor locations placed along the vertical line $x = 0.5$. Specifically, $m = 19$ sensors are uniformly distributed over the interval $y \in [0.2, 1.8]$ (see Fig.~\ref{fig:Schematic_anti}). Observations are assimilated after an initial burn-in period during which the system evolves without measurement updates. Measurements are then discontinued after a prescribed time to examine the evolution of the state in the absence of further observations. Each observation is subject to independent Gaussian noise with zero mean and variance $\sigma_e^2$, where the noise level is chosen as a fixed fraction of the standard deviation of the true displacement signal at the first sensor location. In addition to the observed sensors, the displacement response is monitored at an interior point $(1.75, 1.0)$, which is not included in the observation set. This probe point is used to assess the predictive capability of the filtering framework away from measurement locations.

\subsubsection{Parameter estimation and filtering results}
\label{sec:spe2}
\begin{figure}[!t]
    \centering
	\subfloat[$\vec{x}_p=(0.5,\,0.2)$ \label{fig:forward_anti_1}]
        {\includegraphics[width=0.48\textwidth]{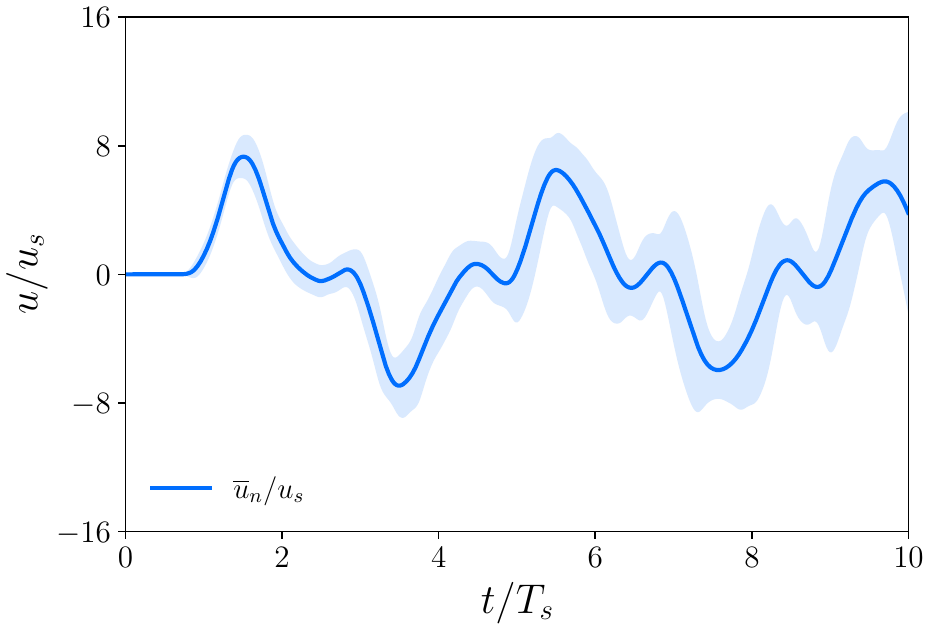}\vspace{-1mm}}\hfill
	\subfloat[$\vec{x}_p=(1.75,\,1.0)$ \label{fig:forward_anti_2}]
        {\includegraphics[width=0.48\textwidth]{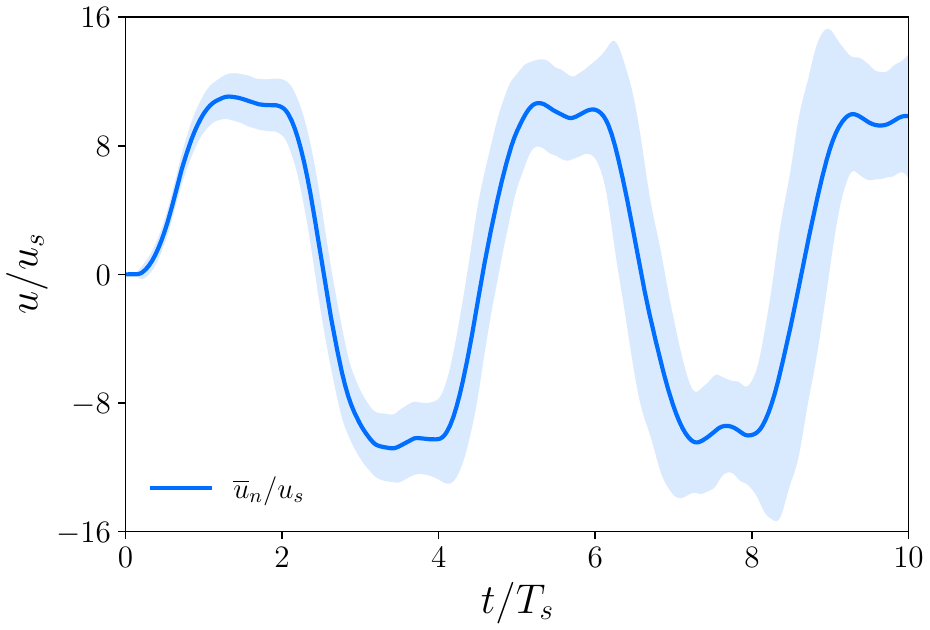}\vspace{-1mm}}
	\caption{Forward predictive displacement response at probe locations $\vec{x}_p$: mean and $95\%$ confidence intervals from the stochastic anti-plane model.}
	\label{fig:forward_anti}
\end{figure}

\begin{figure}[!t]
    \centering
    {\includegraphics[width=0.45\textwidth]{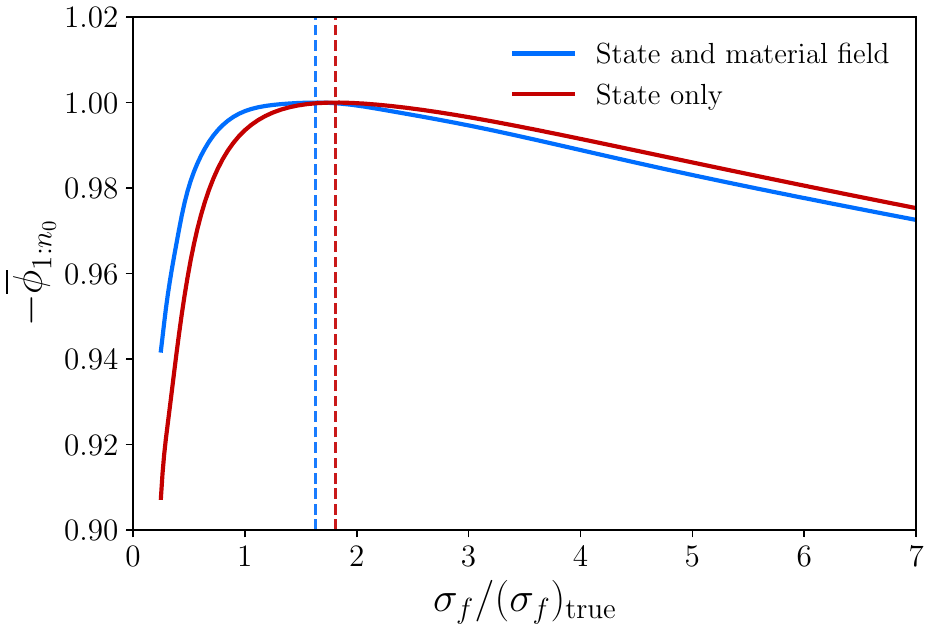}\vspace{-1mm}}
    \caption{Marginal likelihood as a function of the misspecification standard deviation $\sigma_f$.}
    \label{fig:anti_mle}
\end{figure}

The forward predictive response of the stochastic anti-plane model prior to data assimilation is presented in Fig.~\ref{fig:forward_anti}. Figures~\ref{fig:forward_anti_1} and~\ref{fig:forward_anti_2} show the mean displacement and $95\%$ confidence intervals at the observation location $(0.5,0.2)$ and at the point of interest $(1.75,1.0)$, respectively. The predicted mean and variance are noticeably larger at $(1.75,1.0)$, reflecting its proximity to the right boundary where the stochastic traction is applied. In contrast, the observation points lie closer to the clamped boundary at $x=0$ and therefore experience smaller displacement amplitudes and reduced uncertainty. Figure~\ref{fig:anti_mle} presents the estimation of the misspecification standard deviation $\sigma_f$ using the marginal likelihood approach described in Section~\ref{sec:bayes-params}. The likelihood exhibits a well-defined maximum close to the true value of $\sigma_f$, indicating that the misspecification standard deviation is accurately identified from the displacement observations.

\begin{figure}[!t]
    \centering
	\subfloat[$\vec{x}_p=(0.5,\,0.2)$ \label{fig:filter_anti_1}]
        {\includegraphics[width=0.48\textwidth]{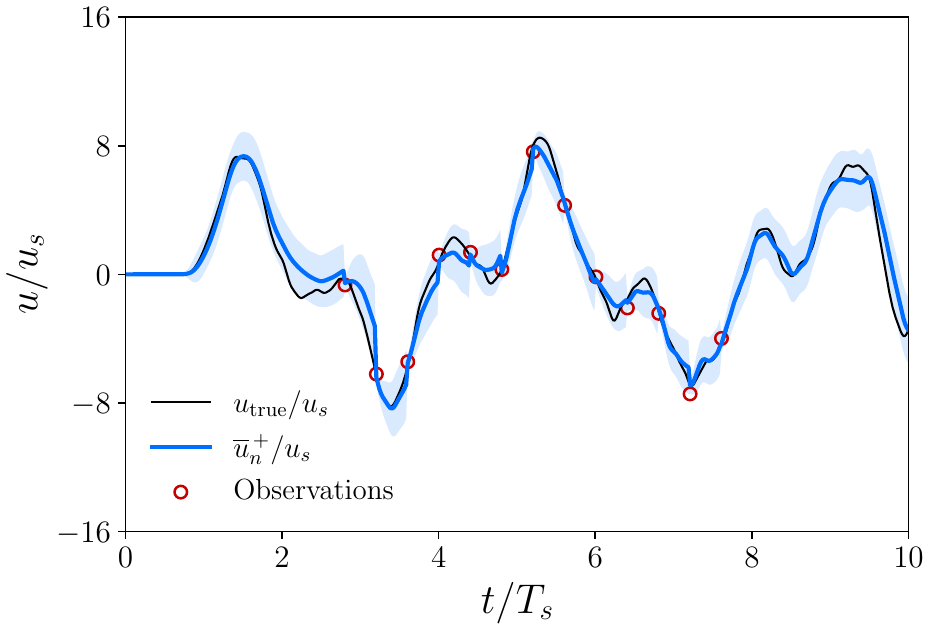}\vspace{-1mm}}\hfill
    \subfloat[$\vec{x}_p=(1.75,\,1.0)$ \label{fig:filter_anti_2}]
        {\includegraphics[width=0.48\textwidth]{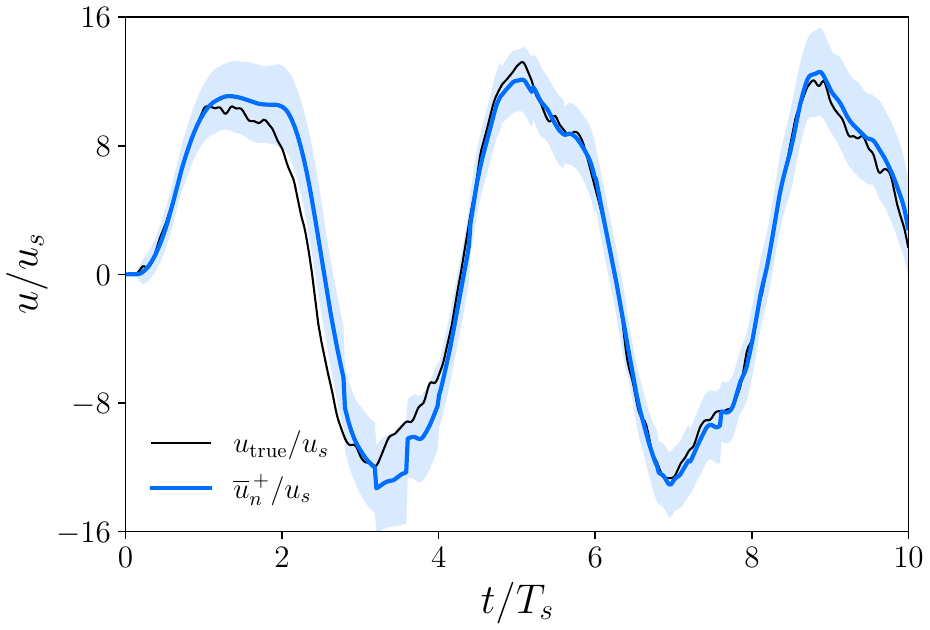}\vspace{-1mm}}
	\caption{Posterior displacement response at probe locations $\vec{x}_p$: mean and $95\%$ confidence intervals after data assimilation.}
	\label{fig:filter_anti}
\end{figure}

Figure~\ref{fig:filter_anti} shows the corresponding posterior displacement estimates after assimilating noisy measurements along the line $x=0.5$, with the misspecification standard deviation $\sigma_f$ taken as the estimated value. In Fig.~\ref{fig:filter_anti_1}, the updated displacement at $(0.5,0.2)$ closely tracks the true displacement, with a marked reduction in uncertainty relative to the forward prediction. More notably, Fig.~\ref{fig:filter_anti_2} demonstrates that the updated displacement at the unobserved interior point $(1.75,1.0)$ also matches the true response with high accuracy. This result highlights the ability of the proposed framework to propagate information from sparsely distributed measurements into regions of the domain that are not directly observed.

\begin{figure}[!t]
    \centering
	\subfloat[Posterior mean $\mu(\vec{x})$ \label{fig:filter_anti_mu_up}]
        {\includegraphics[width=0.48\textwidth]{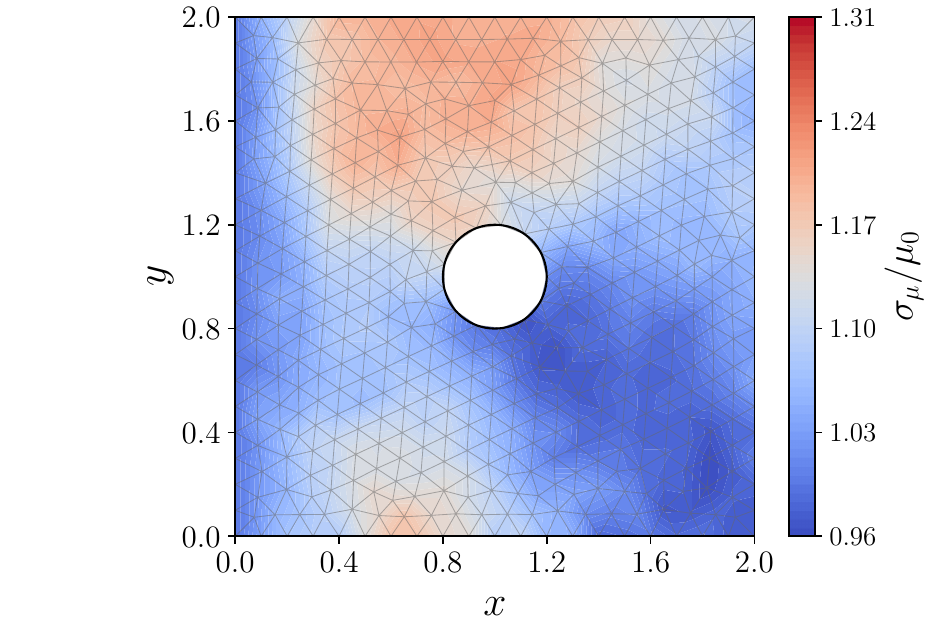}\vspace{-1mm}}
	\subfloat[Posterior standard deviation $\sigma_\mu(\vec{x})$ \label{fig:filter_anti_mu_std}]
        {\includegraphics[width=0.48\textwidth]{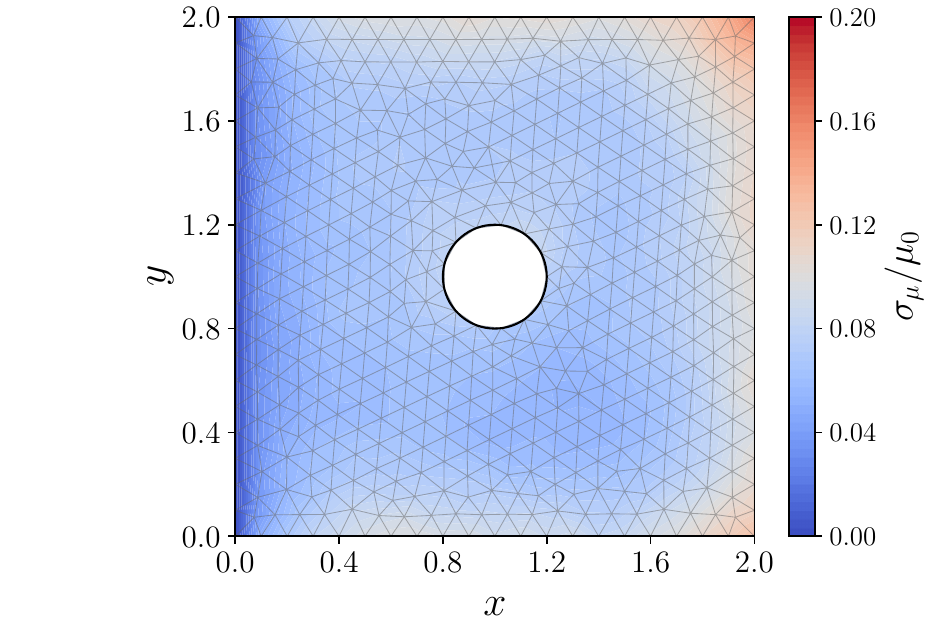}\vspace{-1mm}}
    \vfill
    \subfloat[Shear modulus along $x=0.5$ \label{fig:filter_anti_mu_ob}]
        {\includegraphics[width=0.48\textwidth]{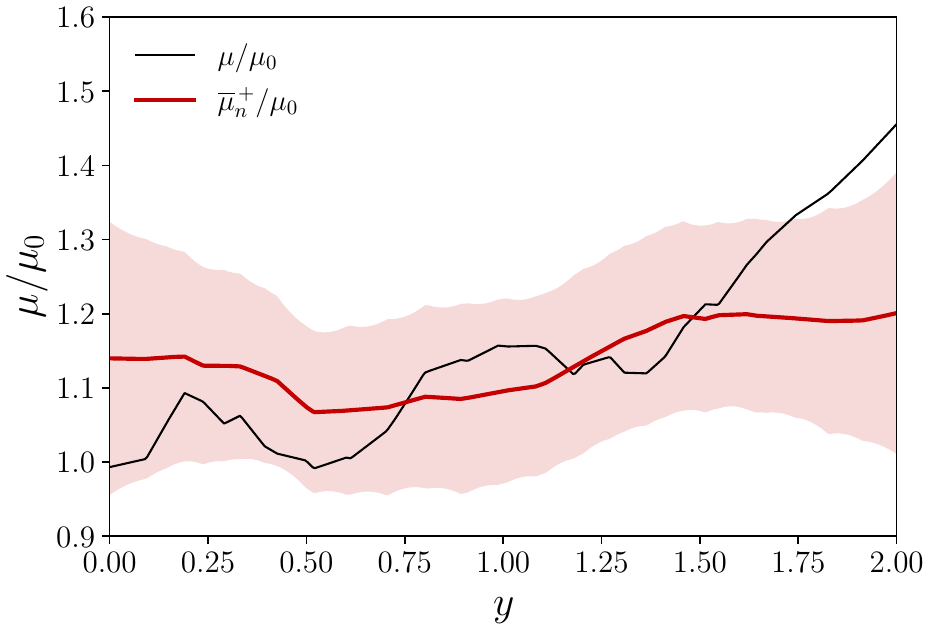}\vspace{-1mm}}
	\caption{Posterior inference of the shear modulus field in the anti-plane problem at the final time.}
	\label{fig:filter_mu}
\end{figure}

The inferred shear modulus field at the final time is shown in Fig.~\ref{fig:filter_mu}. Figures~\ref{fig:filter_anti_mu_up} and~\ref{fig:filter_anti_mu_std} present contour plots of the posterior mean and standard deviation, respectively. Although the prior mean shear modulus is spatially uniform, the updated field exhibits clear spatial structure that closely resembles the true shear modulus distribution shown earlier in Fig.~\ref{fig:shear_field}. In particular, lower values are inferred near the bottom of the domain, while higher values emerge toward the upper region, consistent with the true underlying field. To further quantify the quality of the inference of the material property field, Fig.~\ref{fig:filter_anti_mu_ob} compares the posterior shear modulus along the line $x=0.5$ with the true shear modulus profile along the same line. The posterior mean follows the correct spatial trend, and the true field lies largely within the updated confidence intervals. This agreement demonstrates that the sequential assimilation of displacement data not only improves state estimation but also enables meaningful recovery of spatial material variations starting from a uniform prior.

\begin{figure}[!t]
    \centering
	\subfloat[$\vec{x}_p=(0.5,\,0.2)$ \label{fig:filter_anti_1_so}]
        {\includegraphics[width=0.48\textwidth]{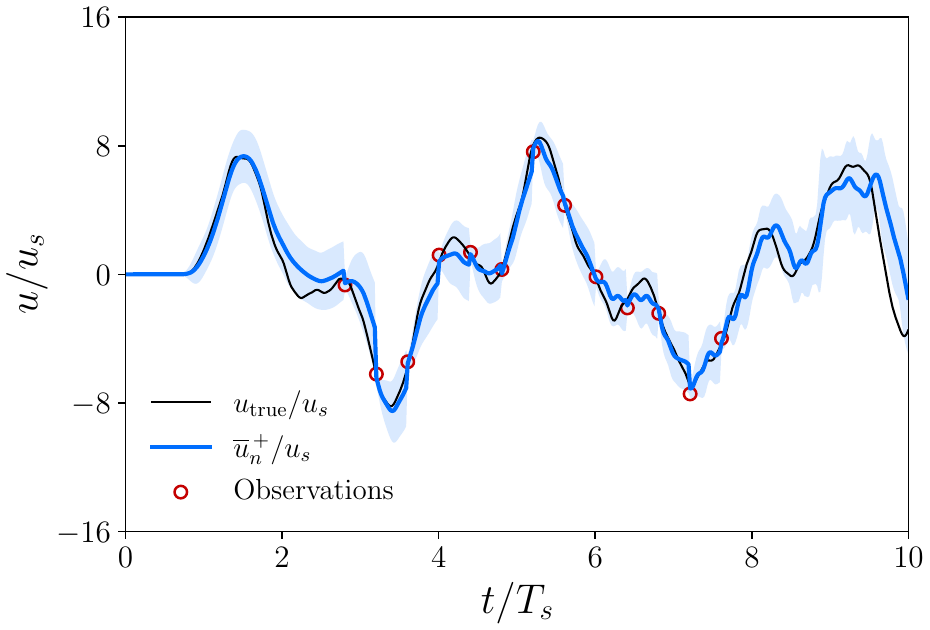}\vspace{-1mm}}\hfill
    \subfloat[$\vec{x}_p=(1.75,\,1.0)$ \label{fig:filter_anti_2_so}]
        {\includegraphics[width=0.48\textwidth]{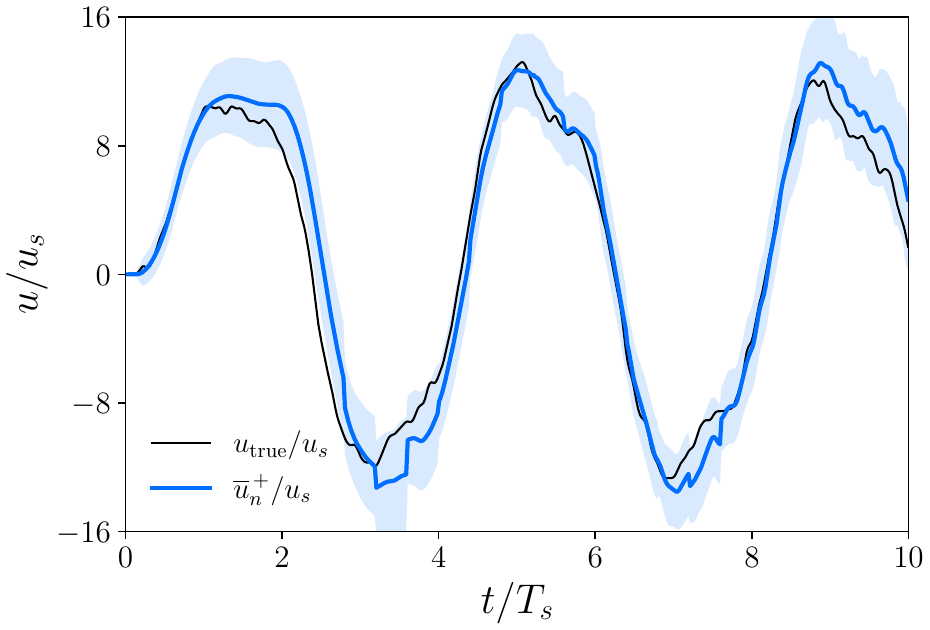}\vspace{-1mm}}
	\caption{Posterior displacement response at probe locations $\vec{x}_p$ without updating the shear modulus field: mean and $95\%$ confidence intervals after data assimilation.}
	\label{fig:filter_anti_so}
\end{figure}

Finally, Fig.~\ref{fig:filter_anti_so} shows the posterior displacement obtained when displacement observations are assimilated while the material field $\kappa$ is not updated, with the misspecification standard deviation $\sigma_f$ taken as the value estimated from the non-augmented formulation. Holding the material field fixed leads to a visibly reduced agreement between the updated and true displacements compared with the joint state and material property filtering results in Fig.~\ref{fig:filter_anti}. This behavior is consistent with the one-dimensional example and indicates that, even in the two-dimensional setting with sparse line measurements, explicitly updating the material field improves the quality of the reconstructed dynamic response. The comparison highlights the role of joint updating of the state and material field in achieving consistent displacement estimates when material variability contributes appreciably to the observed response.

%Overall, the results illustrate that the proposed statFEM framework can accurately infer both dynamic response and spatially varying material properties in an anti-plane elastodynamic setting, even when observations are limited to a sparse set of sensors and the material field is not directly observed.

%% file: conclusions.tex
%
%--------------------------------------------------------------------------------          
\section{Conclusions \label{sec:conclusions}}
%--------------------------------------------------------------------------------
The statistical finite element method (statFEM) provides a probabilistic framework for incorporating observational data and modelling uncertainties into conventional finite element formulations. In this work, we have extended statFEM to sequential inference for transient elastodynamic problems, leading to an augmented Bayesian filtering formulation for joint state and material field estimation. Uncertainties in material properties and forcing are represented as Gaussian spatio-temporal random fields, with model misspecification subsumed into the Brownian forcing component. Their spatial correlation structure is specified through stochastic partial differential equations driven by white noise and discretised consistently with the finite element mesh. To enable tractable filtering in large-scale problems, the nonlinear dependence of the elastodynamic response on the material field is approximated by using a first-order perturbation on the temporally discrete state evolution, yielding a Gaussian approximation to the predictive distribution. This explicit density evolution allows the application of a Bayesian filtering framework in which the augmented prior distribution of the state and material field is propagated through the stochastic dynamics and sequentially updated using noisy displacement observations. The associated hyperparameters governing the random fields are estimated by maximising the marginal likelihood of the observed data.

The numerical examples demonstrate that the proposed framework is able to accurately reconstruct both the displacement field and the underlying material properties in one- and two-dimensional settings. In particular, inferring the material field jointly with the state yields improved posterior estimates compared to formulations in which the material field is held fixed. This improvement is especially pronounced when the material field exhibits strong spatial correlation, where the system response is significantly influenced by coherent material variations. Furthermore, the results show that information from sparse measurements can be effectively propagated to unobserved regions of the domain, enabling accurate prediction of the system response away from sensor locations. More generally, while model misspecification can account for unresolved sources of uncertainty, treating all discrepancies implicitly through forcing leads to less accurate inference. In the present setting of statFEM, explicitly treating the material field as a random field and jointly inferring it together with the state and model misspecification yields more consistent and accurate predictions.

In closing, we outline several possible extensions of the proposed methodology. The first-order perturbation is generally valid for small variances; therefore, more advanced approximations of the evolving state density can be employed, such as polynomial chaos expansion~\cite{narouie2025mechanical,ghanem2003stochastic,sudretPolynomialChaosExpansions2015}. In such settings, application to elastodynamics would require updating the polynomial chaos coefficients~\cite{arnstIdentificationBayesianPosteriors2010} over time~\cite{maiSurrogateModellingStochastic2016,lacourDynamicStochasticFinite2021}. Improved time-marching schemes can also be employed to increase integration order and stability, such as explicit Runge–Kutta methods~\cite{sarkka2019applied} or implicit stochastic Newmark–Beta methods~\cite{bernardStochasticNewmarkScheme2002,royExplorationsFamilyStochastic2005}. Furthermore, employing a spatio-temporal covariance~\cite{sarkkaSpatiotemporalLearningInfiniteDimensional2013} as a physics-informed prior~\cite{weilandFlexibleEfficientProbabilistic2025} for the model misspecification random field may improve state predictions over extended time intervals without additional data~\cite{Tondo2025}. Additionally, reformulating the filtering procedure using SPDE sparse precision matrices that are directly related to the stiffness matrices of elliptic operators~\cite{lindgren2011explicit,chu2021stochastic,ben2024robust} may improve computational efficiency and scalability for large-scale problems. Material uncertainty, typically modelled via scalar random fields, could be more realistically captured using scale-dependent tensor-valued random fields informed by micromechanics ~\cite{torquato2002random,savvas2016determination,jetti2024correlation}. The framework may be extended to scenarios where the constitutive model itself is unknown by coupling the proposed filtering procedure with sparsity-inducing Bayesian model discovery approaches (see, e.g., EUCLID \cite{joshi2022bayesian}), thereby enabling inference of constitutive structure from sparse and dynamic measurements. Finally, approximate Bayesian inference for joint state and parameter estimation may be pursued using variational methods~\cite{SarkaNummenmaa2009} or their stochastic variants~\cite{CourseNair2023}. In this context, a potentially large number of hyperparameters in hierarchical Bayesian filtering may be inferred using scalable stochastic variational algorithms~\cite{HoffmanBlei2013}.

%% file: acknowledgements.tex
\section*{Acknowledgements \label{sec:acknowledgements}}
%--------------------------------------------------------------------------------
%
IK gratefully acknowledges the support of the German Research Foundation (DFG) [Project No. 491258960]. YSJ and FC gratefully acknowledge the support of the Engineering and Physical Sciences Research Council (EPSRC) [Grant Nr. EP/Z003407/1]. AOY acknowledges the support of the Ministry of National Education of the Republic of Turkey.

%% file: appendix.tex
%
%--------------------------------------------------------------------------------          
\section{SPDE Representation of Mat\'ern Random Fields \label{sec:matern}}
%--------------------------------------------------------------------------------
A commonly used choice of operator in the framework described in Section~\ref{sec:forwardgenmat} is a diffusion-type operator, which leads to Mat\'ern random fields. 
Specifically, consider the fractional SPDE
\begin{equation}\label{eq:SPDE}
	\tau \left( \eta^2 - \vec{\nabla}^2 \right)^\beta s(\vec{x}) = \,\xi(\vec{x}),
\end{equation}
where $\vec{\nabla}^2$ is the Laplacian, $\xi(\vec{x})$ is spatial white noise, and the parameters are given by
\begin{equation}
	\eta=\tfrac{\sqrt{2\nu}}{l}, \qquad \beta=\tfrac{\nu}{2}+\tfrac{d}{4}, \qquad 
	{\tau^2}=\frac{\Gamma(\nu)}{{\sigma^2} \,\Gamma(\nu+d/2)(4\pi)^{d/2}(2\nu/l^2)^\nu}.
\end{equation}
Here $d$ is the spatial dimension, $\sigma^2$ the variance, { $\Gamma$ the Gamma function,} $l$ the correlation length and $\nu$ the smoothness parameter. {The exponent can be any value $\beta>1/4$; we select $\beta$ to be an integer by tuning $\nu$ depending on the dimension $d$ to avoid fractional and recursive solution of~\eqref{eq:SPDE}.}

The solution $s(\vec{x})$ is a Gaussian random field with Mat\'ern covariance when the domain is unbounded
\begin{equation}\label{eq:Matern}
	 c_s(\vec{x},\vec{x}^\prime)= \frac{\sigma^2}{2^{\nu-1}\Gamma (\nu)}
	\left(\frac{\sqrt{2\nu}}{l}\|\vec{x}-\vec{x}^\prime\|\right)^{\nu} K_{\nu}\!\left(\frac{\sqrt{2\nu}}{l}\| \vec{x}-\vec{x}^{\prime}\| \right) = \sigma^2 c_{\mathcal{M}}(\vec{x},\vec{x}^\prime),
\end{equation}
where $K_\nu$ is the modified Bessel function of the second kind { of order $\nu$}. The parameter $l$ controls the correlation length, while $\nu$ controls smoothness; in the limit $\nu \to \infty$ the kernel converges to the squared exponential.

For a bounded domain, the discretized form of~\eqref{eq:SPDE} is
\begin{equation}
\mat{L}\, \ary{s}=\mat{M}^{1/2}\ary{\xi},\quad \ary{\xi}\sim\mathcal{N}(\ary{0},\mat{I}),
\end{equation}
{where
\begin{equation}\label{eq:discrete-lap}
\mat{L}=\tau(\eta^2\mat{M}+\mat{K})	
\end{equation}
is a mapping that incorporates the domain discretization and connectivity through the (lumped) Gram matrix $\mat{M}$ and stiffness matrix $\mat{K}$ of the auxiliary problem~\eqref{eq:SPDE}. This yields the covariance 
\begin{equation}
\mat{C}_{\ary{s}}=\mat{L}^{-1}\mat{M}\,\mat{L}^{-\trans}.
\end{equation}	

%
%--------------------------------------------------------------------------------  
\section{Predictive Mean and Covariance \label{sec:predmean}}
%-------------------------------------------------------------------------------
The material parameter $\ary{\kappa}$ is assumed time-invariant, its predictive mean and covariance remains unchanged during the prediction step, i.e.,
\begin{equation}
	\begin{aligned}
\overline{\ary{\kappa}}_{n+1}^{\, -} &=
\mathbb{E}\!\left[\ary{\kappa}| \ary{y}_{1:n}\right] = \overline{\ary{\kappa}}_{n}^{\, +}\\
		\mat{C}_{\kappa|n+1}^{\, -}
&=
\cov(\ary{\kappa}| \ary{y}_{1:n})
=
\mat{C}_{\kappa|n}^{\, +}.
	\end{aligned}
\end{equation}

The predictive mean defined as $\overline{\ary{v}}_{n+1}^{\, -}=\mathbb{E}\!\left[\ary{v}_{n+1}| \ary{y}_{1:n}\right]$
is obtained by taking conditional expectation of the state propagation~\eqref{eq:vvFEM}. Strictly, the mean involves the conditional expectation $\mathbb{E}[\,\mat{A}(\ary{\kappa})\,\ary{v}_n| \ary{y}_{1:n}]$; we adopt a first-order closure and neglect higher-order correlations between $\ary{\kappa}$ and $\ary{v}_n$ under the posterior. Moreover, the process noise increment $\ary{\zeta}_n$ is assumed zero-mean and independent of the observation history $\ary{y}_{1:n}$, so that $\mathbb{E}[\ary{\zeta}_n| \ary{y}_{1:n}]=\mathbb{E}[\ary{\zeta}_n]=\mathbf{0}$.
Setting $\overline{\mat{A}}_n^{\, +} = \mat{A}(\overline{\ary{\kappa}}_n^{\, +})$, we obtain
\begin{equation}\label{eq:appB-mean-expand}
\begin{aligned}
\overline{\ary{v}}_{n+1}^{\, -}
&=
\mathbb{E}\!\left[\left(\mat{A}(\ary{\kappa})\,\ary{v}_n
+ \Delta t\,\mat{B}\,\overline{\ary{f}}_n
+ \ary{\zeta}_n \right) \Big| \ary{y}_{1:n}\right]\\
&=
\mathbb{E}\!\left[\mat{A}(\ary{\kappa})\,\ary{v}_n | \ary{y}_{1:n}\right]
+ \Delta t\,\mat{B}\,\overline{\ary{f}}_n
+ \mathbb{E}\!\left[\ary{\zeta}_n | \ary{y}_{1:n}\right]\\
&\approx
\overline{\mat{A}}_n^{\, +}\,\overline{\ary{v}}_n^{\, +}
+ \Delta t\,\mat{B}\,\overline{\ary{f}}_n.
\end{aligned}
\end{equation}

The predictive covariance defined as
\(
\mat{C}_{n+1}^{\, -}
=
\cov(\ary{v}_{n+1}| \ary{y}_{1:n})
\)
is obtained by evaluating the conditional second moment of the state fluctuation using the linearised update~\eqref{eq:fluc-update}. Adopting a first-order closure consistent with the mean prediction, treating the linearised operators as deterministic, and noting that the process noise increment $\ary{\zeta}_n$ is zero-mean and independent of $(\ary{v}_n,\ary{\kappa})$ and of the observation history $\ary{y}_{1:n}$, we obtain
\begin{equation}\label{eq:appB-cov-expand}
\begin{aligned}
\mat{C}_{n+1}^{\, -}
&=
\mathbb{E}\!\left[
\left(
\overline{\mat{A}}\,(\ary{v}_{n}-\overline{\ary{v}}_n) + \mat{J}_n (\ary{\kappa}-\overline{\ary{\kappa}}) 
+\ary{\zeta}_n
\right)
\left(
\overline{\mat{A}}\,(\ary{v}_{n}-\overline{\ary{v}}_n) + \mat{J}_n (\ary{\kappa}-\overline{\ary{\kappa}}) 
+\ary{\zeta}_n
\right)^\trans
\Big|\,
\ary{y}_{1:n}
\right]\\
&=
\overline{\mat{A}}_n^{\, +}\,\mat{C}_n^{\, +}\,(\overline{\mat{A}}_n^{\, +})^\trans
+\overline{\mat{A}}_n^{\, +}\,\mat{C}_{n,\kappa|n}^{\, +}\,(\mat{J}_n^{\, +})^\trans
+\mat{J}_n^{\, +}\,(\mat{C}_{n,\kappa|n}^{\, +})^\trans\,(\overline{\mat{A}}_n^{\, +})^\trans
+\mat{J}_n^{\, +}\,\mat{C}_{\kappa|n}^{\, +}\,(\mat{J}_n^{\, +})^\trans
+\mat{C}_{\zeta},
\end{aligned}
\end{equation}
where $\mat{J}_n^{\, +}=\mat{J}_n(\overline{\ary{\kappa}}_n^{\, +},\overline{\ary{v}}_n^{\, +})$.

The predictive cross-covariance between the state and the material parameter is obtained by evaluating
\(
\mat{C}_{n+1,\kappa|n}^{\, -}
=
\cov(\ary{v}_{n+1},\ary{\kappa}| \ary{y}_{1:n})
=
\mathbb{E}\!\left[(\ary{v}_{n}-\overline{\ary{v}}_n) (\ary{\kappa}-\overline{\ary{\kappa}})^\trans | \ary{y}_{1:n}\right].
\)
Using the fluctuation update~\eqref{eq:fluc-update}, and invoking the same first-order closure and independence assumptions as in the derivation of $\mat{C}_{n+1}^{\, -}$, we obtain
\begin{equation}
\mat{C}_{n+1,\kappa|n}^{\, -}
=
\overline{\mat{A}}_n^{\, +}\,\mat{C}_{n,\kappa|n}^{\, +}
+
\mat{J}_n^{\, +}\,\mat{C}_{\kappa|n}^{\, +}.
\end{equation}

\section{Updated Mean and Covariance}
\label{sec:app-kalman-update}
%--------------------------------------------------------------------------------

The update step updates the predictive distribution using Bayes’ rule,
\begin{equation}
p(\ary{v}_{n+1},\ary{\kappa}|\ary{y}_{1:n+1})
=
\frac{p(\ary{y}_{n+1}|\ary{v}_{n+1})
\,p(\ary{v}_{n+1},\ary{\kappa}|\ary{y}_{1:n})}
{p(\ary{y}_{n+1}| \ary{y}_{1:n})}.
\label{eq:app-bayes}
\end{equation}
Under the predictive distribution conditioned on $\ary{y}_{1:n}$, the
augmented state vector \(
\left(\ary{v}_{n+1}^\trans\; \ary{\kappa}^\trans\right)^\trans
\)
is Gaussian with mean and covariance given in~\eqref{eq:joint-predictive}.
The observation model is linear,
\begin{equation}
\ary{y}_{n+1} = \mat{H}\ary{v}_{n+1} + \ary{e}_{n+1},
\qquad
\ary{e}_{n+1} \sim \mathcal{N}(\ary{0},\mat{C}_e),
\end{equation}
with $\ary{e}_{n+1}$ independent of $(\ary{v}_{n+1},\ary{\kappa})$.

Consequently, the triplet
\(
(\ary{v}_{n+1},\ary{\kappa},\ary{y}_{n+1})
\)
conditioned on $\ary{y}_{1:n}$ is jointly Gaussian. Its joint density can be written as
\begin{equation}
p\!\left(
\begin{pmatrix}
\ary{v}_{n+1}\\
\ary{\kappa}\\
\ary{y}_{n+1}
\end{pmatrix}
\Bigg|\,
\ary{y}_{1:n}
\right)
=
\mathcal{N}\!\left(
\begin{pmatrix}
\overline{\ary{v}}_{n+1}^{\, -}\\
\overline{\ary{\kappa}}_{n+1}^{\, -}\\
\mat{H}\overline{\ary{v}}_{n+1}^{\, -}
\end{pmatrix},
\begin{pmatrix}
\mat{C}_{n+1}^{\, -} &
\mat{C}_{n+1,\kappa|n}^{\, -} &
\mat{C}_{n+1}^{\, -}\mat{H}^\trans \\[4pt]
(\mat{C}_{n+1,\kappa|n}^{\, -})^\trans &
\mat{C}_{\kappa|n}^{\, -} &
\mat{C}_{\kappa,n+1|n}^{\, -}\mat{H}^\trans \\[4pt]
\mat{H}\mat{C}_{n+1}^{\, -} &
\mat{H}\mat{C}_{n+1,\kappa|n}^{\, -} &
\mat{H}\mat{C}_{n+1}^{\, -}\mat{H}^\trans + \mat{C}_e
\end{pmatrix}
\right).
\end{equation}
The posterior distribution
$\mbox{$p(\ary{v}_{n+1},\ary{\kappa}| \ary{y}_{1:n+1})$}$
can therefore be obtained using the standard conditioning formulas for
block Gaussian distributions. To apply the conditioning formulas, we identify the predictive
moments involving the measurement variable $\ary{y}_{n+1}$. From the observation model, the predictive statistics involving the
measurement are
\begin{align}
\mathbb{E}[\ary{y}_{n+1}|\ary{y}_{1:n}]
&= \mat{H}\overline{\ary{v}}_{n+1}^{\, -},\\
\cov(\ary{y}_{n+1},\ary{y}_{n+1}|\ary{y}_{1:n})
&= \mat{H}\mat{C}_{n+1}^{\, -}\mat{H}^\trans + \mat{C}_e,\\
\cov(\ary{v}_{n+1},\ary{y}_{n+1}|\ary{y}_{1:n})
&= \mat{C}_{n+1}^{\, -}\mat{H}^\trans,\\
\cov(\ary{\kappa},\ary{y}_{n+1}|\ary{y}_{1:n})
&= \mat{C}_{\kappa,n+1|n}^{\, -}\mat{H}^\trans.
\end{align}

Applying the Gaussian conditioning identities~\cite{petersen2008matrix}
to the joint distribution of
$(\ary{v}_{n+1},\ary{\kappa},\ary{y}_{n+1})$
yields the posterior statistics after assimilating the measurement. The updated state mean and covariance are
\begin{align}
\overline{\ary{v}}_{n+1}^{\, +}
&=
\overline{\ary{v}}_{n+1}^{\, -}
+
\mat{C}_{n+1}^{\, -}\mat{H}^\trans
(\mat{H}\mat{C}_{n+1}^{\, -}\mat{H}^\trans + \mat{C}_e)^{-1}
(\ary{y}_{n+1}-\mat{H}\overline{\ary{v}}_{n+1}^{\, -}),\\
\mat{C}_{n+1}^{\, +}
&=
\mat{C}_{n+1}^{\, -}
-
\mat{C}_{n+1}^{\, -}\mat{H}^\trans
(\mat{H}\mat{C}_{n+1}^{\, -}\mat{H}^\trans + \mat{C}_e)^{-1}
\mat{H}\mat{C}_{n+1}^{\, -}.
\end{align}
Similarly, the posterior mean and covariance of the material parameter are
\begin{align}
\overline{\ary{\kappa}}_{n+1}^{\, +}
&=
\overline{\ary{\kappa}}_{n+1}^{\, -}
+
\mat{C}_{\kappa,n+1|n}^{\, -}\mat{H}^\trans
(\mat{H}\mat{C}_{n+1}^{\, -}\mat{H}^\trans + \mat{C}_e)^{-1}
(\ary{y}_{n+1}-\mat{H}\overline{\ary{v}}_{n+1}^{\, -}),\\
\mat{C}_{\kappa|n+1}^{\, +}
&=
\mat{C}_{\kappa|n+1}^{\, -}
-
\mat{C}_{\kappa,n+1|n}^{\, -}\mat{H}^\trans
(\mat{H}\mat{C}_{n+1}^{\, -}\mat{H}^\trans + \mat{C}_e)^{-1}
\mat{H}\mat{C}_{n+1,\kappa|n}^{\, -}.
\end{align}
Finally, the cross-covariance between the state and the material parameter becomes
\begin{equation}
\mat{C}_{n+1,\kappa|n+1}^{\, +}
=
\mat{C}_{n+1,\kappa|n}^{\, -}
-
\mat{C}_{n+1}^{\, -}\mat{H}^\trans
(\mat{H}\mat{C}_{n+1}^{\, -}\mat{H}^\trans + \mat{C}_e)^{-1}
\mat{H}\mat{C}_{n+1,\kappa|n}^{\, -}.
\end{equation}
These expressions correspond to the relations in the update step summarised in
Algorithm~\ref{alg:statFEMDyn}.

\section{Stochastic Verlet Algorithm}\label{sec:app-verlet}
%--------------------------------------------------------------------------------

Several versions of the stochastic Verlet algorithm exist, primarily developed for Langevin dynamics (see e.g.,~\cite{burrageNumericalMethodsSecondOrder2007,mannellaQuasisymplecticIntegratorsStochastic2004a,gronbech-jensenSimpleEffectiveVerlettype2013,melchionnaDesignQuasisymplecticPropagators2007,gronbech-jensenCompleteSetStochastic2020}). Here, we revisit the simple leapfrog method for position-based Verlet time integration~\cite{burrageNumericalMethodsSecondOrder2007,mannellaQuasisymplecticIntegratorsStochastic2004a}.\par 
Without loss of generality, consider the equation of motion of a single-degree-of-freedom system in state space:
\begin{equation}\label{eq:EOM_SDOF_SS}
	\begin{aligned}
		\text{d}u(t) & = \dot{u}(t)\text{d}t \, , \\
		\text{d}\dot{u}(t) &= -\frac{k}{m}u(t)\text{d}t - \frac{\gamma}{m}\dot{u}(t)\text{d}t +\frac{1}{m}\overline{f}(t)\text{d}t + \frac{1}{m}\text{d}\beta(t) \, ,
	\end{aligned}
\end{equation}
where $m$ is the mass, $\gamma$ is the damping coefficient, and $k$ is the stiffness coefficient. The force is resolved in a deterministic (mean) component $\overline{f}$ and a stochastic component through a Brownian motion $\beta(t)$ with diffusion variance $\sigma^2$, such that $\text{d}\beta\text{d}\beta=\sigma_f^2 \text{d}t $.\par 
Integrating the second equation in~\eqref{eq:EOM_SDOF_SS} over a small interval $\Delta t$, between times $t_n$ and $t_{n+1}$, we obtain:
\begin{equation}\label{eq:EOM_SDOF_1}
	\int_{t_{n}}^{t_{n+i}}d\dot{u}(t)  =-\frac{k}{m}\int_{t_n}^{t_{n+1}}u(t)\text{d}t - \frac{\gamma}{m}\int_{t_n}^{t_{n+1}}\dot{u}(t)\text{d}t +\frac{1}{m}\int_{t_n}^{t_{n+1}}\overline{f}(t)\text{d}t + \frac{1}{m}\int_{t_n}^{t_{n+1}}\text{d}\beta(t) \, ,
\end{equation}
which can be approximated as:
\begin{equation}\label{eq:EOM_SDOF_2}
	\dot{u}_{n+1} - \dot{u}_{n} =  - \frac{k \Delta t}{m}  u_{n+1/2} - \frac{\gamma  \Delta t }{m}\dot{u}_{n} + \frac{\Delta t}{m}\overline{f}_n + \frac{1}{m}\Delta\beta_{n} \, 
\end{equation}
where the stochastic integral is approximated considering first-order Taylor expansion, such that $\Delta\beta_n=\mathcal{N}(0,\Delta t\sigma_f^2)$.\par 
The half-step supporting velocity $u_{n+1/2}=u(t_{n+1/2})$ is obtained from~\eqref{eq:EOM_SDOF_SS} as
\begin{equation}\label{eq:EOM_SDOF_3}
	\int_{t_n}^{t_{n+1/2}}du(t) =\int_{t_n}^{t_{n+1/2}} \dot{u}(t)\text{d}t \, ,
\end{equation}
which right-hand side is approximated using left Taylor's first order approximation, resulting in the half-step velocity:
\begin{equation}\label{eq:EOM_SDOF_4}
	\begin{aligned}
		u_{n+1/2} = u_{n} + \frac{\Delta t}{2} \dot{u}_{n}  \, .
	\end{aligned}
\end{equation}\par 
Finally, the displacement at time $t=t_{n+1}$ is completed as:
\begin{equation}\label{eq:EOM_SDOF_5}
	\begin{aligned}
		u_{n+1} = u_{n+1/2} + \frac{\Delta t}{2} \dot{u}_n  \, .
	\end{aligned}
\end{equation}
The integration scheme is commonly conducted by successively applying~\eqref{eq:EOM_SDOF_4}, ~\eqref{eq:EOM_SDOF_2} and~\eqref{eq:EOM_SDOF_5}. These can be reordered in an explicit manner following~\eqref{eq:EOM_SDOF_SS} as:
\begin{equation}
	\begin{pmatrix}
		{u}_{n+1}\\
		\dot{u}_{n+1}
	\end{pmatrix}=
	\begin{pmatrix}\displaystyle
		1 - \frac{1}{2}m^{-1}k\Delta t^2 & \Delta t - \frac{1}{2}m^{-1}\gamma\Delta t ^2 - \frac{1}{4}m^{-1}k \Delta t^3\\
		-m^{-1}k\Delta t & 1-m^{-1}\gamma\Delta t-\frac{1}{2}m^{-1}k \Delta t^2 
	\end{pmatrix}
	\begin{pmatrix}
		u_{n}\\
		\dot{u}_{n}
	\end{pmatrix}
	+
	\begin{pmatrix}
		\frac{1}{2}m^{-1}\Delta t^2\\
		{m^{-1}\Delta t}
	\end{pmatrix}
	\overline{f}_{n}
	+
	\begin{pmatrix}
		\frac{1}{2}m^{-1}\Delta t\\
		{m^{-1}}
	\end{pmatrix}
	\Delta \beta_{n} \, .
\end{equation}
The leapfrog stochastic Verlet for multi-degree-of-freedom systems with mass $\mat{M}$, damping $\mat{D}$ and stiffness $\mat{K}$ matrices is given by:
\begin{equation}
	\ary{v}_{n+1} =
		\begin{pmatrix}
	\mat{I}-\tfrac{\Delta t^2}{2}\mat{M}^{-1}\mat{K} & \Delta t \left( \mat{I}-\tfrac{\Delta t}{2} \mat{M}^{-1}\mat{D}-\tfrac{\Delta t^2}{4}\mat{M}^{-1}\mat{K}\right)\\
	-\Delta t \mat{M}^{-1}\mat{K} & \mat{I}-\Delta t \mat{M}^{-1}\mat{D}-\tfrac{\Delta t^2}{2}\mat{M}^{-1}\mat{K}
\end{pmatrix}\ary{v}_{n}
	+
	\begin{pmatrix}
		\tfrac{\Delta t^2 }{2}\mat{M}^{-1} \\
		\Delta t\mat{M}^{-1}
	\end{pmatrix}
\overline{\ary{f}}_n
	+
	\begin{pmatrix}
		\tfrac{\Delta t }{2}\mat{M}^{-1} \\
		\mat{M}^{-1}
	\end{pmatrix}
	\Delta \ary{\beta}_{n} \, ,
\end{equation}
or put in a compact from:
\begin{equation}
\ary{v}_{n+1}=\mat{A}\,\ary{v}_n 
+ \Delta t\,\mat{B}\,\overline{\ary{f}}_n + \mat{B}\Delta \ary{\beta}_{n} \, .
\end{equation}
The time-marching scheme is quasi-sympletic and retains second-order accuracy for the elastic forces, but is of first order for the damping forces and strong order 1/2 (weak order 1) convergence for the stochastic forces. The stability is guaranteed in the absence of damping, provided that the time step satisfies the usual central-difference condition
\begin{equation}
	\Delta t \leq 2/\omega_{\max} , 
\end{equation}
where $\omega_{\max}$ is the maximum frequency of the discretised system.